\documentclass[12pt,a4paper]{article}
\usepackage[cp1251]{inputenc}

\def\C{\mathbf C}
\def\R{\mathbf R}
\def\N{\mathbf N}
\def\Z{\mathbf Z}
\newcommand{\Cl}{\mathop{\rm Cl}\nolimits}
\newcommand{\Hom}{\mathop{\rm Hom}\nolimits}
\newcommand{\Ann}{\mathop{\rm Ann}\nolimits}
\newcommand{\Supp}{\mathop{\rm Supp}\nolimits}
\newcommand{\Aff}{\mathop{\rm Aff}\nolimits}
\newcommand{\Span}{\mathop{\rm span}\nolimits}
\newcommand{\rank}{\mathop{\rm rank}\nolimits}
\newcommand{\ord}{\mathop{\rm ord}\nolimits}
\newcommand{\adj}{\mathop{\rm adj}\nolimits}
\newcommand{\Id}{\mathop{\rm Id}\nolimits}
\newcommand{\Sol}{\mathop{\rm Sol}\nolimits}

\newtheorem{theorem}{Theorem}[section]
\newtheorem{lemma}[theorem]{Lemma}
\newtheorem{proposition}[theorem]{Proposition}
\newtheorem{corollary}[theorem]{Corollary}
\newtheorem{remark}[theorem]{Remark}
\newcommand{\qedmark}{\ifhmode\unskip\nobreak\hspace{0.5em}\hfill
\else\par\noindent\hfill\fi\rule{1.2mm}{1.2mm}}
\newenvironment{proof}{\par\noindent\textit{Proof:}\enspace
\ignorespaces}{\qedmark\par\medskip}

\makeatletter
\renewcommand{\@seccntformat}[1]{\csname the#1\endcsname.\quad}
\def\@begintheorem#1#2{\trivlist
  \item[\hskip\labelsep{\bfseries #1\ #2.}]\itshape}
\def\@opargbegintheorem#1#2#3{\trivlist
  \item[\hskip\labelsep{\bfseries #1\ #2.\ (#3)}]\itshape}
\let\originalnumberline\numberline
\renewcommand{\numberline}[1]{\originalnumberline{#1.}}
\makeatother

\title{Analytic complexity of solutions of linear\\
PDE systems with constant coefficients\\[4pt]
\large Superpositions of type $n\to k$}
\author{V.~K.~Beloshapka}
\date{}

\begin{document}
\maketitle

\begin{abstract}
For homogeneous systems of linear equations with constant coefficients,
we prove a criterion for finite complexity with respect to superpositions
of type $n\to k$. Namely, each irreducible component of the characteristic
set must lie in an affine subspace of dimension at most $k$.
If this condition fails, the set of entire solutions that admit no
finite representation of type $n\to k$ on any nonempty open subset
is a dense $G_\delta$ subset of the Fr\'echet space of all entire
solutions of the system.
\end{abstract}

\tableofcontents
\newpage

\section*{Introduction}
\addcontentsline{toc}{section}{Introduction}

Analytic complexity of type $n\to k$ measures the depth of a scheme
(the height of a tree) required to represent a function of $n$ variables
by superpositions of functions of $k$ variables, where $k<n$.
If no such representation exists, the complexity is defined to be
infinite. For two variables, the original hierarchy and its
differential-algebraic description were introduced in~\cite{Bel2007};
the multivariable setting is considered in~\cite{Bel2016}.
Further properties of the hierarchy for functions of two variables
were studied in~\cite{Bel2012}. M.~A.~Stepanova~\cite{Step2019}
proved the existence of solutions of the heat equation with infinite
analytic complexity. The present treatment develops the construction
for eliminating internal functions used in~\cite{Step2019} and combines
it with the Baire category theorem and the algebra of systems with
constant coefficients.

For a scalar system, the main algebraic object is the ideal of full
symbols; for several unknown functions, it is a finitely generated
module over the polynomial ring in the frequency variables.
It is essential to use the full symbols: the constant and lower-order
terms of the operators can change the answer. The general theory of
linear systems and its relation to the
Malgrange--Ehrenpreis--Palamodov fundamental principle are presented
in~\cite{AHS,Robertz}. For classical results on the approximation of
solutions by exponential polynomials and on integral representations
of solutions, see~\cite{Mal1956,Ehr1970,Pal1970}. A description of primary
ideals by means of Noetherian differential operators is given
in~\cite{CHS2021}.

The argument rests on the affine geometry of the irreducible components
of the spectral support. Multiplicities affect the value of the finite
upper bound on complexity. The exact class indices also depend on the
arity of the composition tree and on the convention for level zero.
These three issues are treated separately in the statements.
The main results concern holomorphic solutions of homogeneous linear
systems with constant coefficients.

\section{Scalar systems and superpositions of type $n\to k$}
\label{sc:core}

In this section we study systems with one unknown function
(the scalar case).

\subsection{Superposition classes and the solution space}
\label{sc:definitions}

Let $n\geq2$, $1\leq k<n$, $X=(x_1,\ldots,x_n)$, and
$D_i=\partial/\partial x_i$. All functions and representations are
considered locally in a complex domain. For $k\geq2$, let
$\Cl^{n,k}_0$ denote the class of holomorphic functions depending on
at most $k$ original coordinates, and set
\[
 \Cl^{n,k}_{\ell+1}
 =\{G(A_1,\ldots,A_k): A_j\in\Cl^{n,k}_{\ell},\quad
 G\ \mbox{holomorphic}\}.
\]
We call the parameter $\ell$ the \textit{depth}, and the functions used
to define the composition the \textit{internal} functions.
The classes are nested and contain the constants.

For $k=1$, addition must be included among the allowed operations.
We adopt the following convention for depth. The class $\Cl^{n,1}_0$
consists of functions of one original coordinate,
\[
 \Cl^{n,1}_1
 =\left\{g\left(\sum_{i=1}^n a_i(x_i)\right):
 g,a_1,\ldots,a_n\ \mbox{holomorphic}\right\},
\]
and, for $\ell\geq1$, we set
\[
 \Cl^{n,1}_{\ell+1}
 =\{g(A+B): A,B\in\Cl^{n,1}_{\ell},\quad
 g\ \mbox{holomorphic}\}.
\]
For $n=2$, this is the usual definition of binary analytic complexity.
Other conventions for the number of summands at the first level may
change the numerical value of the depth, but not the finiteness of
the representation. Throughout what follows, the term
``depth'' refers precisely to these classes.

For a holomorphic function $u$ on a connected domain $\Omega$, let
$N_{n,k}(u)$ denote the smallest depth of a representation on some
nonempty open subdomain of $\Omega$. If there is no such representation,
set $N_{n,k}(u)=\infty$. Thus, infinite complexity means that no finite
representation exists on any nonempty open subset.
Local solutions are considered on arbitrarily small connected
neighborhoods; depth estimates concern a representation at a generic
point of such a neighborhood. We do not require the same representation
to exist at a prescribed exceptional point. Whenever representations
are added or multiplied below, a common domain of definition is chosen.
No additional assumptions are used to compare exact depths on different
open subsets.

Finite representations are closed under addition, multiplication, and
composition with functions of at most $k$ arguments. For $k\geq2$,
addition and multiplication are admissible internal functions.
For $k=1$, we use the identity
\[
 AB=\frac{(A+B)^2-A^2-B^2}{2}.
\]
Linear forms also have finite representations. Hence a linear change
of the original coordinates preserves finiteness of complexity,
although it need not preserve its exact value.

Let
\[
 R_n=\C[\xi_1,\ldots,\xi_n],\quad
 I=(P_1,\ldots,P_s)\subset R_n,\quad
 V(I)=\{\xi\in\C^n:P_j(\xi)=0,\ 1\leq j\leq s\},
\]
so that $I$ is an ideal and $V(I) \subseteq \mathbf{C}^n$ is its
common zero set.

We consider the scalar homogeneous system
\begin{equation}
 P_j(D)u=0,\qquad 1\leq j\leq s.
 \label{sc:system}
\end{equation}
The ideals $I=(0)$ and $I=R_n$ are also allowed. Let $\mathcal E(I)$
denote the space of entire solutions on $\C^n$, equipped with the
topology of uniform convergence on compact sets. It is a closed
linear subspace of the Fr\'echet space of entire functions:
differentiation is continuous in this topology.
Clearly, $\exp(\xi \cdot x) \in \mathcal E(I)$ if and only if
$\xi \in V(I)$.
The affine algebraic variety $V(I)$ is called the {\it characteristic}
or {\it spectral} set of~(\ref{sc:system}).

For an irreducible algebraic set $W\subset\C^n$, set
\[
 d_{\rm aff}(W)=\dim\Aff(W).
\]
This is the dimension of the smallest affine linear subspace containing
$W$. In general, $\dim W\leq d_{\rm aff}(W)$, and the inequality
may be strict.

\begin{theorem}[Scalar criterion]
\label{sc:criterion}
For system~(\ref{sc:system}), the following conditions are equivalent.
\begin{enumerate}
\item All local holomorphic solutions have finite complexity of type
$n\to k$.
\item There is a common finite upper bound on the complexity of all
local holomorphic solutions.
\item Every irreducible component $W$ of $V(I)$ satisfies
\[
 d_{\rm aff}(W)\leq k.
\]
\end{enumerate}
The last condition is equivalent to the inclusion of $V(I)$ in a
finite union of affine linear subspaces of dimension at most $k$.

If this condition fails, the set of entire solutions of infinite
complexity contains a dense $G_\delta$ subset of $\mathcal E(I)$.
Moreover, finite sums of characteristic exponentials have arbitrarily
large finite complexity.
\end{theorem}

The necessity of the geometric condition and the final assertions of
the theorem are proved in this section. Sufficiency, including the
common upper bound, will be proved in Proposition~\ref{pm:positive}.

\subsection{Characteristic exponentials and independent coordinates}
\label{sc:coordinates}

Fix an irreducible component $W\subset V(I)$ such that
\[
 d=d_{\rm aff}(W)>k.
\]
Let $\mathcal I(W)$ denote its ideal, and let
\[
 A=\C[\xi_1,\ldots,\xi_n]/\mathcal I(W)
\]
be its coordinate ring. This ring is an integral domain. Set
\[
 \mathcal T_W=
 \left\{\sum_{\nu=1}^N c_\nu e^{\langle\lambda^{(\nu)},X\rangle}:
 N<\infty,\quad c_\nu\in\C,\quad \lambda^{(\nu)}\in W\right\}.
\]
Here $\langle\lambda,X\rangle=\sum_i\lambda_i x_i$ is the bilinear,
not the Hermitian, pairing. Every function in $\mathcal T_W$ is an
entire solution of the original system.

For $f\in A$ and $w\in\mathcal T_W$, define
\[
 Z_f[w](X)=\sum_{\nu=1}^N
 c_\nu e^{\langle\lambda^{(\nu)},X\rangle} f(\lambda^{(\nu)}).
\]
If $f$ is represented by a polynomial, then $Z_f[w]=f(D)w$.
The definition is independent of the representative of the class $f$.
We have the identities
\begin{equation}
 D_i Z_f[w]=Z_{\xi_i f}[w].
 \label{sc:shift-indices}
\end{equation}

\begin{lemma}[Independence of coordinates]
\label{sc:independence}
If $f_1,\ldots,f_N\in A$ are linearly independent over $\C$, then,
for every fixed point $X_0\in\C^n$, the map
\[
 \mathcal T_W\longrightarrow\C^N,\qquad
 w\longmapsto (Z_{f_1}[w](X_0),\ldots,Z_{f_N}[w](X_0))
\]
is surjective. A finite set of frequencies suffices, and this set can
be chosen independently of $X_0$.
\end{lemma}

\proof
The functionals $f\mapsto f(\lambda)$, $\lambda\in W$, span the dual
space of $\Span_\C(f_1,\ldots,f_N)$. Otherwise, a nonzero linear
combination of the $f_i$ would vanish on all of $W$ and hence would
be zero in $A$. Therefore there exist
$\lambda^{(1)},\ldots,\lambda^{(N)}\in W$ such that the matrix
\[
 \bigl(f_i(\lambda^{(j)})\bigr)_{i,j=1}^N
\]
is nonsingular. Passing to the values $Z_{f_i}[w](X_0)$ multiplies
its columns by the nonzero numbers
$e^{\langle\lambda^{(j)},X_0\rangle}$. An appropriate choice of
the coefficients $c_j$ therefore yields any vector in $\C^N$.
\endproof

\subsection{A vector field annihilating the arguments at a vertex}
\label{sc:fields}

Consider the collection of linear functions that are constant on $W$.
\[
 K_W=\{a\in\C^n:a\cdot\xi\ \mbox{is constant on }W\}.
\]
This is a linear space. If $L_W$ is the direction space of the affine
subspace $\Aff(W)$, then $K_W=L_W^\perp$ with respect to the
bilinear pairing. Consequently,
\[
 \dim K_W=n-d.
\]

\begin{lemma}[Choosing a vector field outside $K_W$]
\label{sc:field}
Let $s_1,\ldots,s_p$ be holomorphic functions on a connected domain
$U\subset\C^n$, where $0\leq p\leq k$. There exists a nonzero
vector field
\[
 \Delta=\sum_{i=1}^n a_i(X)D_i
\]
such that $\Delta s_j=0$ for every $j$ and $a(X)\notin K_W$ at
generic points. The coefficients of the field can be chosen as
polynomial expressions with integer coefficients in the first
derivatives of the $s_j$. For fixed $n,p$, there is a finite list
of such expressions containing a field with the stated properties.
\end{lemma}

\proof
Let $r$ be the maximal rank of the matrix $J=(D_i s_j)$; its rank
equals $r$ on a nonempty open set. If $r=0$, the list of coordinate
fields $D_i$ suffices. At least one of their directions lies outside
$K_W$, since $d>k\geq1$.

Suppose $r>0$. Choose $r$ rows and $r$ columns forming a square
submatrix $M$ whose determinant is not identically zero. Let
$i_1,\ldots,i_r$ be the indices of the selected columns. For every
free column with index $j$, let $b_j$ be its restriction to the
selected rows, and set
\begin{equation}
 \Delta_j=(\det M)D_j-
 \sum_{\nu=1}^r\bigl((\adj M)b_j\bigr)_\nu D_{i_\nu},
 \label{sc:minor-field}
\end{equation}
where $(\adj M)$ is the adjugate matrix of $M$. The identity
$M\adj M=(\det M)\Id$ shows that the field annihilates the selected
functions $s_i$. On the set of constant rank, the remaining rows are
meromorphic linear combinations of the selected rows, and hence
$\Delta_j s_i=0$ for all arguments. These equalities extend as
analytic identities.

At a point where $\det M\ne0$, the coefficient vectors of the
fields~(\ref{sc:minor-field}) form a basis of $\ker J$. However,
\[
 \dim\ker J=n-r\geq n-k>n-d=\dim K_W.
\]
Thus at least one vector of this basis lies outside $K_W$.
Choose such a vector at one nonexceptional point. Its membership in
$K_W$ is expressed by a finite system of analytic linear equations;
therefore, outside a proper analytic subset, the selected vector
remains outside $K_W$. All choices of rank, rows, minors, and free
columns form a finite list. The formula involves no division by
$\det M$ and uses only integer coefficients.
\endproof

For $n=3$, $k=2$, and rank two, one may take
\[
 \Delta=(\nabla s_1\times\nabla s_2)\cdot(D_1,D_2,D_3).
\]
For smaller rank, the corresponding minors of smaller order are used.
For a vertex depending directly on at most $k$ coordinates, the lemma
selects a suitable coordinate vector field corresponding to a missing
coordinate.

\subsection{Isolating the terms of highest pole order}
\label{sc:valuation}

Let $Y$ be the normalization of the projective closure of $W$.
The coordinates $\xi_i$ are rational functions on $Y$.
For every $a\notin K_W$, the function $q=a\cdot\xi$ is nonconstant
and therefore has a pole along some irreducible divisor $E\subset Y$.
Indeed, on a normal variety, a rational function that is regular at
every point of codimension one is regular everywhere. A global regular
function on an irreducible projective variety is constant.

For a nonzero rational function $f$, set
\[
 h_E(f)=-\ord_E(f),\qquad h_E(0)=-\infty.
\]
Then
\begin{equation}
 h_E(fg)=h_E(f)+h_E(g),\qquad
 h_E(f+g)\leq\max(h_E(f),h_E(g)).
 \label{sc:valuation-rules}
\end{equation}
The orders may be negative.

If $a=a(X)$ consists of the coefficients of a field from
Lemma~\ref{sc:field}, choose a nonexceptional point $X_0$ at which
$a(X_0)\notin K_W$, and a pole divisor of the function
$a(X_0)\cdot\xi$. At generic points,
\begin{equation}
 h_E(q(X,\xi))=s>0,\qquad q(X,\xi)=\sum_i a_i(X)\xi_i,
 \label{sc:q-pole}
\end{equation}
where $s$ is independent of $X$. To verify this, it suffices to
consider the finite-dimensional space $\Span_\C(\xi_1,\ldots,\xi_n)$
with its filtration by $h_E$. Vanishing of the projection onto any
of its graded pieces is expressed by a finite system of linear
conditions on the coefficients $a_i(X)$. The highest piece with a
projection that is not identically zero gives the constant generic
order $s$.

\begin{lemma}[Absence of rational first integrals]
\label{sc:no-first-integrals}
Let
\[
 \Delta=\sum_i a_i(X)D_i
\]
be a meromorphic vector field on a connected domain $U$, with
$a(X)\notin K_W$ generically. Let $R[w]$ be a rational expression in
finitely many derivatives of $w\in\mathcal T_W$, with fixed
meromorphic coefficients in $X$. The denominator is assumed not to
vanish identically on this family. If $R$ depends essentially on the
test function, then $\Delta R\not\equiv0$ on $\mathcal T_W$.

The same assertion holds for every translated family
$u_*+\mathcal T_W$, where $u_*$ is a fixed holomorphic function on $U$.
\end{lemma}

\proof
First consider $\mathcal T_W$. All derivatives involved are indexed
by a finite-dimensional subspace $S\subset A$. Choose a divisor $E$
for which~(\ref{sc:q-pole}) holds. In $S$, choose a constant complex
basis adapted to the filtration
\[
 S\cap\{f:h_E(f)\leq j\}.
\]
This filtration has only finitely many nonzero successive quotients.
Its intersection as $j\to-\infty$ is zero, and finite-dimensionality
of $S$ therefore allows us to choose a basis successively at all
levels. The initial forms of the basis elements at each level are
linearly independent over $\C$.

By Lemma~\ref{sc:independence}, the corresponding values
$z_i=Z_{f_i}[w]$ are independent coordinates. Write $R$ as an element
of the field $\mathcal M(U)(z_1,\ldots,z_N)$, where $\mathcal M(U)$
is the field of meromorphic functions on $U$, cancel common factors,
and discard the inessential variables. Let $h$ be the largest order
among the indices of the remaining variables. Denote the highest
nonzero terms by $f_1,\ldots,f_b$, so that $h_E(f_i)=h$.
The choice of basis and~(\ref{sc:valuation-rules}) imply, for every
fixed nonexceptional $X$,
\[
 h_E\left(\sum_{i=1}^b\lambda_i qf_i\right)=h+s
 \quad\mbox{if}\quad (\lambda_1,\ldots,\lambda_b)\ne0.
\]
Hence $qf_1,\ldots,qf_b$ are independent modulo all functions of order
less than $h+s$. In particular, they are independent modulo the
previous essential indices and the products $qf_j$ from lower levels.

The chain rule and~(\ref{sc:shift-indices}) give
\begin{equation}
 \Delta R=\Delta_{\rm coeff}R+\sum_i R_{z_i}Z_{qf_i}[w]
 =\sum_{i=1}^b R_{z_i}Z_{qf_i}[w]+T.
 \label{sc:highest-group}
\end{equation}
Here $\Delta_{\rm coeff}$ denotes differentiation of the coefficients
with the $z_i$ held fixed; in this formula, $Z_{qf_i}$ means
$\sum_j a_j(X)Z_{\xi_j f_i}$. All frequency indices of the coordinates
on which $T$ and the coefficients $R_{z_i}$ depend have order less
than $h+s$. We have $R_{z_i}\ne0$ for at least one $i\leq b$.

Fix a nonexceptional $X$ at which this partial derivative remains
a nonzero rational expression. Applied to the combined finite set
of indices, Lemma~\ref{sc:independence} allows the new group of
highest-order terms to vary while the old coordinates are held fixed.
Expression~(\ref{sc:highest-group}), which is linear in this new
group, therefore cannot vanish identically.

For the translated family, write $v=u_*+w$. The substitution
\[
 D^\alpha v=D^\alpha u_*+D^\alpha w
\]
changes only the meromorphic coefficients of the rational expression.
Further differentiation also places derivatives of $u_*$ in the
coefficient part. The free coordinates always refer to $w$ and still
satisfy~(\ref{sc:shift-indices}). There is no need here to apply the
class $f\in A$ as a differential operator to $u_*$ itself.
The preceding argument remains valid without change.
\endproof

\subsection{Elimination of internal functions}
\label{sc:elimination}

\begin{proposition}[An additional autonomous relation]
\label{sc:elimination-proposition}
Let $u_*$ be a holomorphic function of finite complexity of type
$n\to k$ on a nonempty connected domain $U$. Then there exists a
differential polynomial
\begin{equation}
 F\in\Z[u_\alpha:\alpha\in\Z_{\geq0}^n],\qquad F[u_*]=0,
 \label{sc:integer-relation}
\end{equation}
that is nontrivial on the translated test family $u_*+\mathcal T_W$.
Each such polynomial involves only finitely many derivatives.
\end{proposition}

\proof
Fix a finite representation of $u_*$. If necessary, duplicate
identical subexpressions to obtain a tree with no identified vertices.
Absorb all numerical parameters into the internal functions.
Number the vertices so that each ancestor precedes its descendants.

A vertex has the form $g(s_1,\ldots,s_p)$, where $p\leq k$.
For $k=1$, the single argument may be a sum of child functions;
at the first level, a sum of $n$ functions of one variable is allowed.
This changes neither the chain rule nor the fact that the argument
depends only on descendants. Functions at the leaves are regarded
as vertices whose arguments are the selected original coordinates.
At this stage, even a polynomial internal function is retained as
a separate vertex.

For each vertex, choose a field from Lemma~\ref{sc:field}.
Its coefficients depend only on descendants. It annihilates all
values of partial derivatives of the current internal function:
\[
 \Delta\bigl((\partial^\beta g)(s_1,\ldots,s_p)\bigr)=0
 \quad\mbox{for all }\beta.
\]
This equality also holds when the arguments have deficient rank,
including the case of constant arguments.

Set $v=u_*+w$, $w\in\mathcal T_W$. The field associated with the
root gives the first relation
\[
 H[v]=\Delta_{\rm root}v,\qquad H[u_*]=0.
\]
It is nontrivial on the family $u_*+\mathcal T_W$: its variation
with respect to $w$ is $Z_q[w]$, where $q=a(X)\cdot\xi$ is
nonconstant on $W$. The root function has already been eliminated;
only its descendants remain in the coefficients.

Suppose that some vertices have been eliminated and that an
expression $H[v]$ has been obtained which is nontrivial on
$u_*+\mathcal T_W$, polynomial in the derivatives of $v$ and in
the values of partial derivatives of the remaining internal
functions, and satisfies $H[u_*]=0$. Collect its formal monomials
with respect to the current vertex:
\begin{equation}
 H[v]=\sum_{j=1}^N\alpha_j[v]F_j(s(X)),\qquad
 \Delta F_j(s(X))=0.
 \label{sc:collected}
\end{equation}
The coefficients $\alpha_j$ do not involve the current internal
function. Nontriviality of $H$ ensures the existence of a coefficient
$\alpha_1$ that is nontrivial on the translated family.
Consider the rational expression $R=H/\alpha_1$ on this family.

If $R$ depends essentially on the free jet coordinates,
Lemma~\ref{sc:no-first-integrals} gives $\Delta R\not\equiv0$.
Define a new polynomial expression by
\begin{equation}
 H_1[v]=\sum_{j=2}^N
 \bigl(\alpha_1[v]\Delta\alpha_j[v]-\alpha_j[v]\Delta\alpha_1[v]\bigr)
 F_j(s(X)).
 \label{sc:wronskian-step}
\end{equation}
Fix the representation of $u_*$ and vary only the test function
$v=u_*+w$, where $w\in\mathcal T_W$. Since
$\Delta(F_j(s(X)))=0$ for all $j$, the Leibniz rule
applied to~(\ref{sc:collected}) gives
\[
 H_1[v]=\alpha_1[v]\,\Delta(H[v])-H[v]\,\Delta(\alpha_1[v]).
\]
The term containing $F_1$ cancels on the right-hand side;
this explains why the sum in~(\ref{sc:wronskian-step}) starts
at $j=2$.

At $v=u_*$, we have $H[u_*]\equiv0$, and hence also
$\Delta(H[u_*])\equiv0$. It follows that $H_1[u_*]\equiv0$.
On the other hand, in the field of rational functions in the free
jet coordinates, the last identity for $H_1[v]$ takes the form
\[
 H_1=\alpha_1^2\,\Delta(H/\alpha_1).
\]
Both factors on the right-hand side are nonzero: the first by the
choice of $\alpha_1$, and the second by
Lemma~\ref{sc:no-first-integrals}. Thus $H_1$ does not vanish
identically on the test family.

If $R$ is independent of the free jet coordinates, then on the
entire translated family
\[
 H[v]=C(X)\alpha_1[v],\qquad C\not\equiv0.
\]
This is an equality of polynomial expressions in the free coordinates,
so it can be evaluated at $w=0$, even if the original fraction
$H/\alpha_1$ was undefined there. We obtain $\alpha_1[u_*]=0$.
At the same time, $\alpha_1$ remains nontrivial on the family and
no longer involves the current vertex. This step completes its
elimination.

Repeat the procedure. In the first branch, the number of monomials
of the current vertex strictly decreases; when only one monomial
remains, the ratio $H/\alpha_1$ is already independent of $v$.
Thus the vertex is eliminated after finitely many steps.
Differentiating $\alpha_j$ and the coefficients of $\Delta$
introduces derivatives only of the remaining vertices and their
descendants. An eliminated vertex never reappears.

All expressions retained in the algorithm are obtained by addition,
subtraction, multiplication, differentiation, selection of a minor,
and selection of a coefficient $\alpha_1$.
Formula~(\ref{sc:wronskian-step}) is applied to the fixed
representation and the chosen rank of its arguments; the identities
$\Delta s_j=0$ do not require solving any linear equations that
introduce new numerical coefficients. The chain rule produces only
integer combinatorial coefficients. The divisor, frequencies, and
basis are chosen solely to verify nontriviality and do not alter
the output polynomial.

Once all internal functions have been eliminated, only derivatives
of $v$ remain in the expression. The original coordinates, which
serve as arguments at the leaves, do not remain explicitly either:
the coefficients in the algorithm used only their derivatives,
which are zero or one. We have obtained the autonomous
polynomial~(\ref{sc:integer-relation}) with integer coefficients,
retaining both required properties.
\endproof

\begin{lemma}[A finite list at bounded depth]
\label{sc:finite-list}
For fixed $n,k,\ell$, there is a finite list of autonomous
differential polynomials with integer coefficients having the
following property. For every function $u_*$ of depth at most
$\ell$, at least one polynomial in this list vanishes on $u_*$
and is nontrivial on $u_*+\mathcal T_W$.
\end{lemma}

\proof
Under the chosen definition of depth, the number of vertices of the
tree is bounded in terms of $n,k,\ell$. There are only finitely many
vertex types, ways of connecting them, and sets of original coordinates
at the leaves. It therefore suffices to consider one fixed formal tree.

For each of its vertices, the field formulas from
Lemma~\ref{sc:field} form a finite list: all ranks and all choices of
rows and columns are allowed. For each formal expression already
obtained, the list of monomials of the current vertex is finite.
At step~(\ref{sc:wronskian-step}), this list only decreases:
the coefficients $\alpha_j$, their derivatives, and the coefficients
of the field do not involve the current vertex. Thus, starting with
$N$ monomials, at most $N-1$ such steps are possible before the final
selection of a coefficient. At each step there are only finitely many
choices of the field and the coefficient.

Consequently, for a fixed input expression, there are only finitely
many possible outputs after the elimination of one vertex.
Although derivatives of descendants and the number of their monomials
may increase, each output is a specific finite polynomial expression.
Induction on the finite number of vertices yields finitely many final
polynomials. If coefficients or monomials vanish for a special
representation, this merely selects a sublist of the existing choices.
The numbers specifying the internal functions themselves are not
substituted into this formal procedure.
Proposition~\ref{sc:elimination-proposition} guarantees that, for
every actual representation, the finite list contains a polynomial
with the required properties.
\endproof

\begin{corollary}[A polynomial for one level]
\label{sc:level-polynomial}
For every $\ell$, there is an autonomous differential polynomial
$F_\ell$ with integer coefficients, nontrivial on $\mathcal T_W$,
that vanishes on every function $v\in\mathcal T_W$ admitting a local
representation of depth at most $\ell$.
\end{corollary}

\proof
If $u_*\in\mathcal T_W$, then $u_*+\mathcal T_W=\mathcal T_W$.
From the finite list in Lemma~\ref{sc:finite-list}, remove the
polynomials that vanish identically on $\mathcal T_W$ and multiply
the remaining ones. For every function of the stated depth,
at least one factor vanishes on that function. The product is nontrivial: after passing to
independent frequency coordinates for the finite union of all indices
involved, it is a product of nonzero polynomials over $\C$.
A local identity for an entire function extends to all of $\C^n$.
\endproof

\begin{corollary}[A relation on the compact-open closure]
\label{sc:normal-closure}
Suppose that an entire function $u$ belongs to the closure of
$\mathcal T_W$ in the topology of uniform convergence on compact
subsets of $\C^n$ and has finite complexity on some nonempty open
set. Then there is an autonomous differential polynomial $F$ with
integer coefficients such that $F[u]\equiv0$ and $F$ is nontrivial
on $\mathcal T_W$.
\end{corollary}

\proof
Proposition~\ref{sc:elimination-proposition} gives a polynomial that
is nontrivial on $u+\mathcal T_W$. If it vanished identically on
$\mathcal T_W$, continuity of derivatives under locally uniform convergence
would extend this identity to the entire closure of $\mathcal T_W$.
The latter is a linear space containing $u+\mathcal T_W$, which
gives a contradiction. The equality $F[u]=0$ extends from the domain
of the representation to all of $\C^n$.
\endproof

\subsection{The Baire category theorem in the space of entire solutions}
\label{sc:baire}

Suppose that the geometric condition of Theorem~\ref{sc:criterion}
fails, and let $W$ be the component chosen above. For every entire
solution $u_*$ of finite complexity,
Proposition~\ref{sc:elimination-proposition} gives a polynomial
$F$ that is nontrivial on $\mathcal E(I)$, since
\[
 u_*+\mathcal T_W\subset\mathcal E(I).
\]
The local equality $F[u_*]=0$ is an identity on all of $\C^n$.

Enumerate all autonomous differential polynomials with integer
coefficients that are nontrivial on $\mathcal E(I)$ as
$F_1,F_2,\ldots$. There are at most countably many of them. Set
\[
 \mathcal Z_j=\{u\in\mathcal E(I):F_j[u](0)=0\}.
\]
Each $\mathcal Z_j$ is closed, since a finite set of derivatives
at a fixed point depends continuously on an entire function in the
topology of uniform convergence on compact sets.

Nontriviality of $F_j$ means that there exist a solution $v$ and
a point $X_0$ with $F_j[v](X_0)\ne0$. The system has constant
coefficients and $F_j$ is autonomous. Translating the argument,
we may therefore assume that $F_j[v](0)\ne0$.
For arbitrary $u\in\mathcal E(I)$, the expression
\[
 F_j[(1-t)u+tv](0)
\]
is a polynomial in $t$ that is nonzero at $t=1$.
It cannot vanish for all sufficiently small $t$. Moreover,
$(1-t)u+tv\to u$ in $\mathcal E(I)$ as $t\to0$.
Hence $\mathcal Z_j$ has empty interior.

By the Baire category theorem, the set
\begin{equation}
 \mathcal R=\mathcal E(I)\setminus\bigcup_{j\geq1}\mathcal Z_j
 \label{sc:residual-set}
\end{equation}
is a dense $G_\delta$ subset of $\mathcal E(I)$.
Every entire solution admitting a finite representation on a nonempty
open set belongs to one of the sets $\mathcal Z_j$.
Thus all functions in $\mathcal R$ have infinite complexity.
In particular, this proves the necessity of the geometric condition
in Theorem~\ref{sc:criterion}.

Finally, for every $\ell$, the polynomial $F_\ell$ from
Corollary~\ref{sc:level-polynomial} is nontrivial on $\mathcal T_W$.
By autonomy and translation invariance of this family, there exists
$v\in\mathcal T_W$ with $F_\ell[v](0)\ne0$.
Such a function has no representation of depth at most $\ell$ on any
open set. At the same time, it is a finite sum of exponentials of
linear forms, and therefore has finite complexity.
Thus the finite complexities of characteristic exponential sums
are unbounded.

\section{Analytic additivity of solution spaces}
\label{pm:analytic}

Throughout this section,
\[
 R=\C[\xi_1,\ldots,\xi_n],\qquad
 \xi_jf=\frac{\partial f}{\partial z_j}.
\]
We write $\mathcal O(\Omega)$ for the space of holomorphic functions on an
open domain $\Omega\subset\C^n$, and $\mathcal O_0$ for the space of
holomorphic germs at the origin. All modules below are considered over
the polynomial ring $R$.

We use the following classical result on systems with constant
coefficients: if $G\subset\R^N$ is an open convex domain, then
$C^\infty(G,\C)$ is an injective module over
$\C[\partial_{x_1},\ldots,\partial_{x_N}]$; see~\cite{AHS}, Section~2,
and the sources on the fundamental principle cited there. Injectivity
means that, for every inclusion of modules $B\subset C$, every module
homomorphism $B\to C^\infty(G,\C)$ extends to $C$. It is this result,
rather than an assertion about the convergence of a formal solution,
that is used below.

\begin{lemma}\label{pm:injectivity}
If $\Omega\subset\C^n$ is convex and open, then $\mathcal O(\Omega)$
is an injective $R$-module. In particular, this holds for $\mathcal O(\C^n)$.

For holomorphic germs, the following assertion, sufficient for our
applications, holds: given an inclusion of finitely generated $R$-modules
$B\subset C$, every homomorphism $B\to\mathcal O_0$ extends to $C$.
\end{lemma}

\proof
Set $z_j=x_j+iy_j$ and consider the ring
\[
 R'=\C[\xi_1,\ldots,\xi_n,\eta_1,\ldots,\eta_n],
\]
acting on smooth functions by
\[
 \xi_j=\partial_{z_j}=\frac12(\partial_{x_j}-i\partial_{y_j}),
 \qquad
 \eta_j=\partial_{\bar z_j}=\frac12(\partial_{x_j}+i\partial_{y_j}).
\]
This is an invertible complex linear change of the $2n$ real
derivations. Hence the quoted result applies to
$C^\infty(\Omega,\C)$ as an $R'$-module.

Let $B\subset C$ be an arbitrary inclusion of $R$-modules, and let
$\varphi:B\to\mathcal O(\Omega)$ be an $R$-linear map. Make
$B$ and $C$ into $R'$-modules by setting $\eta_jB=\eta_jC=0$.
Since the values of $\varphi$ are holomorphic, the map
\[
 \varphi:B\longrightarrow C^\infty(\Omega,\C)
\]
is $R'$-linear. By injectivity, it extends to an
$R'$-linear map $\psi:C\to C^\infty(\Omega,\C)$. For every
$c\in C$, we have
\[
 \partial_{\bar z_j}\psi(c)=\psi(\eta_jc)=0
 \qquad(j=1,\ldots,n).
\]
A smooth function satisfying all the Cauchy--Riemann equations is
holomorphic. Thus the image of $\psi$ is contained in $\mathcal O(\Omega)$.
This proves injectivity for arbitrary modules, not merely for cyclic
modules of the form $R/J$.

Now consider a homomorphism $\varphi:B\to\mathcal O_0$, where
$B\subset C$ are finitely generated. Choose a finite set of generators
$b_1,\ldots,b_q$ of $B$. Their images are represented by holomorphic
functions $f_1,\ldots,f_q$ on a common sufficiently small ball. Since
$R$ is Noetherian, the module of relations among $b_1,\ldots,b_q$
has a finite set of generators. The corresponding differential relations
among $f_1,\ldots,f_q$ hold as equalities of germs. Shrinking the ball
finitely many times, we ensure that all these relations hold on a common
convex ball $\Omega$. We thus obtain an $R$-linear map
$B\to\mathcal O(\Omega)$ representing $\varphi$. Extend it to
$C$ by the assertion already proved and pass to germs. The resulting
map $C\to\mathcal O_0$ extends $\varphi$.
\endproof

\begin{remark}\label{pm:germs-remark}
The local assertion of the lemma already gives the exactness of the
functor $\Hom_R(\_,\mathcal O_0)$ on finitely generated modules.
This suffices for all subsequent local applications. Baire's theorem
is applied below in spaces of holomorphic functions on a fixed domain,
in particular in the Fr\'echet space of entire solutions; no topology
on the set of all germs is required for this purpose.
\end{remark}

For an ideal $J\subset R$, set
\[
 \Sol_\Omega(J)=\{u\in\mathcal O(\Omega):P(D)u=0
 \mbox{ for all }P\in J\}.
\]
The space of germs $S(J)$ and the space of entire solutions
$\mathcal E(J)$ are defined similarly.

\begin{lemma}\label{pm:additivity}
For a finite collection of ideals $J_1,\ldots,J_s\subset R$ and an
open convex domain $\Omega\subset\C^n$, we have
\[
 \Sol_\Omega\left(\bigcap_{\alpha=1}^sJ_\alpha\right)
 =\sum_{\alpha=1}^s\Sol_\Omega(J_\alpha).
\]
The same equality holds for spaces of holomorphic germs and for
spaces of entire solutions.
\end{lemma}

\proof
The diagonal map gives an embedding
\[
 R/\bigcap_\alpha J_\alpha\longrightarrow
 \bigoplus_{\alpha=1}^s R/J_\alpha.
\]
By Lemma~\ref{pm:injectivity}, applying
$\Hom_R(\_,\mathcal O(\Omega))$ gives a surjection
\[
 \bigoplus_{\alpha=1}^s\Hom_R(R/J_\alpha,\mathcal O(\Omega))
 \longrightarrow
 \Hom_R\left(R/\bigcap_\alpha J_\alpha,\mathcal O(\Omega)\right).
\]
A homomorphism from $R/J$ is determined by the image of the class
of the unit element; the admissible images of the unit element form precisely
$\Sol_\Omega(J)$. Under these identifications, the resulting surjection
is $(u_1,\ldots,u_s)\mapsto\sum_\alpha u_\alpha$. This proves the
required inclusion; the reverse inclusion follows directly from
$\bigcap_\alpha J_\alpha\subset J_\alpha$. For entire functions, take
$\Omega=\C^n$; for germs, use the local part of
Lemma~\ref{pm:injectivity}.
\endproof

\section{Sufficiency of the affine spectral condition}
\label{pm:positive-section}

Let $1\le k<n$. Complexity is understood in the sense of the
superpositions $n\to k$ introduced above: holomorphic functions
of at most $k$ arguments are allowed. When $k=1$, addition must
also be allowed. Multiplication can be expressed in terms of
addition and functions of one variable:
\[
 fg=\frac14((f+g)^2-(f-g)^2).
\]
Consequently, the class of functions of finite complexity is closed
under finite sums and products, as well as under the substitution
of linear forms for the original coordinates.

\begin{lemma}\label{pm:thick-plane}
Let $L\subset\C^n$ be an affine subspace of dimension $d$,
$J=I(L)$, and $m\ge1$. After a linear change of independent
variables, there exist coordinates $(s,t)\in\C^d\times\C^c$,
$c=n-d$, and a vector $b\in\C^c$ in which all local solutions
of $J^m(D)u=0$ have the form
\begin{equation}\label{pm:normal-form}
 u(s,t)=e^{b\cdot t}\sum_{|\beta|<m}t^\beta f_\beta(s).
\end{equation}
The functions $f_\beta$ are arbitrary holomorphic functions. For
entire solutions, they are entire. Their number is
$ {m+c-1\choose c}$; when $d=0$, they are constants.
\end{lemma}

\proof
Choose linear frequency coordinates $(\sigma,\tau)$ in which $L$
is defined by the equations $\tau_i=b_i$, $i=1,\ldots,c$.
Let $(s,t)$ be the dual independent coordinates, so that the pairing
of the frequency and the independent variable takes the form
$\sigma\cdot s+\tau\cdot t$. Then
\[
 J=(\tau_1-b_1,\ldots,\tau_c-b_c).
\]
The ideal $J^m$ is generated by the monomials
$\prod_i(\tau_i-b_i)^{\beta_i}$, $|\beta|=m$.
The substitution $u=e^{b\cdot t}v$ transforms the corresponding
equations into
\[
 \partial_t^\beta v=0\qquad(|\beta|=m).
\]
All derivatives with respect to $t$ of total order at least $m$
also vanish. The Taylor series in $t$ near any point $(s_0,t_0)$
therefore terminates. Expanding the powers $(t-t_0)^\beta$, we obtain
\[
 v(s,t)=\sum_{|\beta|<m}t^\beta f_\beta(s).
\]
Its coefficients are holomorphic in $s$. If $u$ is entire, then
$f_\beta(s)=\partial_t^\beta v(s,0)/\beta!$ is also entire.
Substitution back into the equations verifies all of them directly.
The number of multi-indices $\beta\in\N_0^c$ with $|\beta|<m$
is $ {m+c-1\choose c}$.
\endproof

\begin{proposition}\label{pm:positive}
Let $I\subset R$, and suppose that every irreducible component
of $V(I)$ is contained in an affine subspace of dimension at most
$k$. Then the complexities of all local holomorphic solutions of
the system $I(D)u=0$ are bounded by a single finite number depending
only on $I,n,k$. Entire solutions admit a representation as a finite
sum of expressions of the form~(\ref{pm:normal-form}) with entire
functions $f_\beta$.
\end{proposition}

\proof
If $I=R$, the only solution is zero. Otherwise, take a primary
decomposition
\[
 I=\bigcap_{\alpha=1}^sQ_\alpha,\qquad \mathcal P_\alpha=\sqrt{Q_\alpha}.
\]
The irreducible set $V(\mathcal P_\alpha)$ is contained in $V(I)$.
By assumption, the latter is covered by finitely many affine
subspaces of dimension at most $k$. An irreducible set contained
in a finite union of closed sets is contained in one of them.
Hence, for every $\alpha$, we can choose an affine subspace
$L_\alpha$ such that
\[
 V(\mathcal P_\alpha)\subset L_\alpha,\qquad d_\alpha=\dim L_\alpha\le k.
\]
This argument also applies to embedded primary components.

Set $J_\alpha=I(L_\alpha)$. Then
$J_\alpha\subset \mathcal P_\alpha$. There is an integer
$m_\alpha\ge1$ such that
$\mathcal P_\alpha^{m_\alpha}\subset Q_\alpha$.
Indeed, choose generators $g_1,\ldots,g_r$ of
$\mathcal P_\alpha$ and exponents $a_i\ge1$ with
$g_i^{a_i}\in Q_\alpha$. For $m_\alpha=1+\sum_i(a_i-1)$,
every monomial of degree $m_\alpha$ in $g_1,\ldots,g_r$
is divisible by at least one $g_i^{a_i}$ and thus belongs
to $Q_\alpha$. Therefore,
\[
 J_\alpha^{m_\alpha}\subset Q_\alpha.
\]

By Lemma~\ref{pm:additivity}, every local solution can be represented
as $u=\sum_\alpha u_\alpha$, where $Q_\alpha(D)u_\alpha=0$.
Each $u_\alpha$ also satisfies
$J_\alpha^{m_\alpha}(D)u_\alpha=0$ and therefore has the normal
form of Lemma~\ref{pm:thick-plane}. In this form, the arbitrary
functions depend on $d_\alpha\le k$ linear forms, while the
remaining factors are exponentials of linear forms and monomials.
All these expressions are constructed by admissible superpositions.
Their number and structure depend on the ideal and the chosen
primary decomposition, but not on the solution. For entire solutions,
additivity is also global, so the coefficient functions in the
normal forms can be chosen to be entire.

For an explicit rough bound, set
\[
 c_\alpha=n-d_\alpha,\qquad
 T=\sum_{\alpha=1}^s{m_\alpha+c_\alpha-1\choose c_\alpha},
 \qquad b_n=\lceil\log_2n\rceil.
\]
A linear form in the original coordinates is constructed using at
most $b_n$ binary levels of addition. A function of one variable
applied to such a form, and a function of at most $k$ such forms,
are constructed using at most $b_n+1$ levels. A monomial in at
most $n$ linear forms is constructed by raising each form to the
corresponding power and using balanced multiplication. If $f$
and $g$ admit representations of depths $a$ and $b$ on the same
open set, then $fg$ admits a representation there of depth at
most $\max(a,b)+2$. Using this estimate, we obtain the bound
$3b_n+1$ for a monomial. Further multiplication by an exponential
and by a function $f_\beta$ adds at most four levels. Summing
at most $T$ terms adds at most $\lceil\log_2T\rceil$ levels.
Thus, under the adopted numbering of the classes, the following
bound suffices:
\begin{equation}\label{pm:uniform-bound}
 N_{n,k}(u)\le3\lceil\log_2n\rceil+5+\lceil\log_2T\rceil.
\end{equation}
This is certainly not optimal, but it is manifestly independent of
$u$. With another equivalent numbering convention, only the
technical depth bound changes, not uniform finiteness.
\endproof

The nonreduced structure of the ideal enters into the exponents
$m_\alpha$, the number of terms, and the numerical depth bound.
The condition for uniform finiteness itself is determined by the
reduced spectrum. The proof of the proposition does not require
integration of arbitrary functions of finite complexity with respect
to one variable: it uses analytic additivity and a finite Taylor
series in the transverse directions.

\section{Matrix systems with several unknowns}
\label{pm:matrix-section}

Let
\begin{equation}\label{pm:matrix-system}
 A(D)U=0,\qquad
 A\in\C[\xi_1,\ldots,\xi_n]^{p\times m},\qquad
 U=(u_1,\ldots,u_m)^t.
\end{equation}
Write
\[
 N=\mathop{\rm im}A^t\subset R^m,\qquad M=R^m/N.
\]
Here $A^t$ denotes the transpose matrix. The classes of the standard
basis vectors in $M$ are denoted by $e_1,\ldots,e_m$.
The space of entire vector solutions is denoted by
$\mathcal E(A)=\Sol_{\C^n}(A)$. Finite complexity of a vector
solution on an open set means finite complexity of all its components
on that set. A more precise form of the assertion about infinite
complexity will be obtained below: one fixed component admits no
finite representation on any nonempty open set.

\begin{lemma}\label{pm:module-description}
There is a natural identification
\[
 \Sol_\Omega(A)=\Hom_R(M,\mathcal O(\Omega)).
\]
It also holds for germs and for entire functions. The spectral
support of the module is given by
\begin{equation}\label{pm:support}
 V_A=\Supp M=V(\Ann_RM)=
 \{\lambda\in\C^n:\rank A(\lambda)<m\}.
\end{equation}
If $p\ge m$, the last set is defined by all minors of order
$m$ of the matrix $A$. If $p<m$, it equals $\C^n$.
\end{lemma}

\proof
The map $R^m\to\mathcal O(\Omega)$ sending the standard basis
vectors to the corresponding $u_j$ factors through $M$ if and
only if all rows of $A(D)U$ vanish. This proves the first
identification.

For the maximal ideal
$\mathcal M_\lambda=(\xi_1-\lambda_1,\ldots,\xi_n-\lambda_n)$,
the fiber $M/\mathcal M_\lambda M$ is isomorphic to the cokernel
of $A(\lambda)^t$. By Nakayama's lemma~\cite{AM1969}, the
localization of the finitely generated module $M$ at
$\mathcal M_\lambda$ is nonzero if and only if this fiber is
nonzero. The latter condition is equivalent to $\rank A(\lambda)<m$.
The support of a finitely generated module coincides with the zero
set of its annihilator. The rank condition is equivalent to the
vanishing of all maximal minors when $p\ge m$; when $p<m$,
the rank is always less than $m$.
\endproof

Set
\[
 I_j=\Ann_R(e_j)=\{P\in R:Pe_j=0\mbox{ in }M\}.
\]
Since $e_1,\ldots,e_m$ generate $M$, we have
\begin{equation}\label{pm:coordinate-annihilators}
 \Ann_RM=\bigcap_{j=1}^mI_j,\qquad
 V_A=\bigcup_{j=1}^mV(I_j).
\end{equation}

\begin{lemma}\label{pm:projection}
For every $j$ and every open convex domain $\Omega\subset\C^n$,
the projection
\[
 \pi_j:\Sol_\Omega(A)\longrightarrow\Sol_\Omega(I_j),
 \qquad U\longmapsto u_j,
\]
is surjective. The same holds for germs and for entire solutions.
\end{lemma}

\proof
The cyclic submodule $Re_j\subset M$ is isomorphic to $R/I_j$.
Any solution $f\in\Sol_\Omega(I_j)$ defines a homomorphism
$Re_j\to\mathcal O(\Omega)$ sending $e_j$ to $f$. By
Lemma~\ref{pm:injectivity}, it extends to $M$. The corresponding
vector solution has $u_j=f$. For germs, the local part of the
same lemma is used.
\endproof

\begin{theorem}\label{pm:matrix-criterion}
Let $1\le k<n$. For the matrix system~(\ref{pm:matrix-system}),
the following conditions are equivalent:
\begin{enumerate}
\item all local holomorphic solutions have components of finite
complexity $n\to k$;
\item the complexities of all components of all local solutions are
bounded by a single number depending only on $A,n,k$;
\item every irreducible component $W$ of $V_A$ satisfies
\[
 \dim\Aff(W)\le k.
\]
\end{enumerate}
If the last condition fails, there is an index $j$ for which
the set of entire vector solutions with $N_{n,k}(u_j)=\infty$
contains a dense $G_\delta$ subset of $\mathcal E(A)$.
Moreover, such a solution can be chosen to be of exponential type.
\end{theorem}

\proof
Suppose that condition~3 holds and $I=\Ann_RM$.
Every component of any solution of the system satisfies
$I(D)u_j=0$. By Proposition~\ref{pm:positive}, the complexities
of these components have a common finite upper bound. This proves
$3\Rightarrow2$. The implication $2\Rightarrow1$ is immediate.

Let $W$ be an irreducible component of $V_A$ with
$\dim\Aff(W)>k$. Over the function field $\C(W)$, the matrix
$A$ has rank less than $m$. Hence the space
\[
 M\otimes_R\C(W)=\C(W)^m/\mathop{\rm im}A^t
\]
is nonzero. Choose $j$ so that the image of $e_j$ in this
space is nonzero. Then $I_j\subset I(W)$: if a polynomial
in $I_j$ did not vanish identically on $W$, it would be
invertible in $\C(W)$, and its annihilation of the image
of $e_j$ would contradict that image being nonzero.
Thus $W\subset V(I_j)$. Since $V(I_j)\subset V_A$, and
$W$ is an irreducible component of $V_A$, the set $W$
is also a component of $V(I_j)$.

Apply the scalar criterion, Theorem~\ref{sc:criterion}, to the
ideal $I_j$. We obtain an entire function
$f\in\mathcal E(I_j)$ of infinite complexity.
By Lemma~\ref{pm:projection}, it lifts to $U\in\mathcal E(A)$
with $u_j=f$. This contradicts condition~1 and proves the
equivalence.

We show that the solutions for which the fixed component $u_j$
has finite complexity form a set of the first Baire category
in $\mathcal E(A)$. For sets whose complements are of the first
Baire category, we use the term \emph{residual}, in accordance
with the English terminology. We show that the solutions for
which the fixed component $u_j$ has infinite complexity form
a residual subset of $\mathcal E(A)$. Indeed, the spaces
$\mathcal E(A)$ and $\mathcal E(I_j)$ are closed linear
subspaces of $\mathcal O(\C^n)^m$ and $\mathcal O(\C^n)$,
respectively, equipped with the topology of uniform convergence
on compact sets. They are Fr\'echet spaces.
By Lemma~\ref{pm:projection}, the map
$\pi_j:\mathcal E(A)\to\mathcal E(I_j)$ is continuous,
linear, and surjective. By the open mapping theorem, it is open.

The scalar criterion gives a dense $G_\delta$ set,
$G=\bigcap_{\nu\ge1}G_\nu\subset\mathcal E(I_j)$, consisting
of functions of infinite complexity; the sets $G_\nu$ can
be chosen to be open and dense. Each $\pi_j^{-1}(G_\nu)$
is open and dense: the image of any nonempty open set in
$\mathcal E(A)$ is open and therefore intersects $G_\nu$.
Consequently,
\[
 \pi_j^{-1}(G)=\bigcap_{\nu\ge1}\pi_j^{-1}(G_\nu)
\]
is dense and is of type $G_\delta$. For all its elements,
the same component $u_j$ has infinite complexity.

The existence of a solution of exponential type with this property
will be proved directly in Proposition~\ref{pm:exp-type}.
The index $j$ can be kept the same: the nonzero image of
$e_j$ in $M\otimes_R\C(W)$ admits a linear functional
whose value on $e_j$ is~1, giving a polarization with
$j$th component equal to~1.
\endproof

\begin{corollary}\label{pm:local-statement}
If the condition of Theorem~\ref{pm:matrix-criterion} fails,
then in a neighborhood of every point of $\C^n$ there exists
a holomorphic vector solution with one fixed component admitting
no finite representation on any nonempty open subset of that
neighborhood.

Moreover, for every nonempty open convex domain
$\Omega\subset\C^n$, the set of such solutions contains a dense
$G_\delta$ subset of $\Sol_\Omega(A)$ in the topology of
uniform convergence on compact sets. The component index $j$
can be chosen independently of $\Omega$.
\end{corollary}

\proof
For the first assertion, it suffices to restrict the entire
solution from the theorem to the required neighborhood.
For the second, retain $W$ and $j$ chosen in its proof.
Let $\mathcal H_\Omega$ be the countable set of autonomous
differential polynomials with integer coefficients that are
nontrivial on $\Sol_\Omega(I_j)$. For
$F\in\mathcal H_\Omega$, set
\[
 Z_F=\{f\in\Sol_\Omega(I_j):F[f]\equiv0\mbox{ on }\Omega\}.
\]
This set is closed. It has empty interior: choose
$v\in\Sol_\Omega(I_j)$ and $z_0\in\Omega$ such that
$F[v](z_0)\ne0$. For any $f\in Z_F$, the expression
$F[(1-t)f+tv](z_0)$ is a polynomial in $t$ that is nonzero
at $t=1$. It is therefore nonzero for arbitrarily small
$t$, while $(1-t)f+tv\to f$ in the specified topology.

If $f\in\Sol_\Omega(I_j)$ has a finite representation on
some nonempty open subset,
Proposition~\ref{sc:elimination-proposition} gives an autonomous
polynomial $F$ with integer coefficients such that $F[f]=0$
on that subset and which is nontrivial on $f+\mathcal T_W$.
By the identity theorem, $F[f]=0$ throughout the connected
domain $\Omega$. Moreover,
$f+\mathcal T_W\subset\Sol_\Omega(I_j)$, since
$I_j\subset I(W)$. Hence $F\in\mathcal H_\Omega$.
Baire's theorem shows that the complement of
$\bigcup_{F\in\mathcal H_\Omega}Z_F$ is dense and is of type
$G_\delta$; none of its elements has a finite representation
on any nonempty open subset.

Finally, the projection
$\pi_j:\Sol_\Omega(A)\to\Sol_\Omega(I_j)$ is surjective
by Lemma~\ref{pm:projection} and open by the open mapping
theorem for Fr\'echet spaces. The preimage of the set just
constructed is the required dense $G_\delta$ set.
\endproof

The matrix formulation of the qualitative criterion follows
from the scalar theorem and the injectivity of the space of
analytic solutions. The ideals $I_j$ have a precise analytic
meaning: they describe all admissible individual components
of a solution, since the projection in
Lemma~\ref{pm:projection} is surjective.

\section{Polarizations and solutions of exponential type}
\label{pm:exp-section}

We prove the last part of the theorem with control of growth.
An entire vector solution is said to be of exponential type
if there exist $C,T>0$ such that
\[
 \|U(z)\|\le Ce^{T\|z\|}\qquad(z\in\C^n).
\]
Any fixed norms may be used in this definition. Below we use
Euclidean norms. Consider the coefficient space
$X=\prod_{r\ge1}[2^{-r-1},2^{-r}]$ with the topology of
coordinatewise convergence. To each point $c\in X$, we associate
a solution $U_c$.

\begin{proposition}\label{pm:exp-type}
Let $W\subset V_A$ be an irreducible component with
$\dim\Aff(W)>k$, and suppose that the image of $e_j$ in
$M\otimes_R\C(W)$ is nonzero. Then there exists an entire
solution of $A(D)U=0$ of exponential type whose $j$th component
has infinite complexity $n\to k$. More precisely, this property
holds on a dense $G_\delta$ subset of $X$.
\end{proposition}

\proof
\emph{1. A normalized polarization.}
Choose a linear functional over $\C(W)$ on
$M\otimes_R\C(W)$ taking the value~1 on the image of $e_j$.
Let $b_1,\ldots,b_m\in\C(W)$ be its values on the standard
generators. Then
\begin{equation}\label{pm:polarization}
 A(\lambda)b(\lambda)=0,\qquad b_j(\lambda)=1.
\end{equation}
These equalities are understood as identities of rational
functions on $W$. There is a nonempty Zariski open subset
of the smooth locus of $W$ on which all $b_i$ are regular.
Choose in it a small closed coordinate polydisc $K$
containing a nonempty relatively open part of $W$.
Both frequencies and polarizations are bounded on $K$:
\[
 \|\lambda\|\le T,\qquad\|b(\lambda)\|\le B\qquad(\lambda\in K),
\]
where we may assume $T\ge1$ and $B\ge1$.

\emph{2. A sequence of frequencies and its tails.}
Choose pairwise distinct points $\lambda_r\in K$, $r\ge1$,
forming a dense subset of $K$, such that every nonempty
relatively open subset of $K$ contains infinitely many of
the chosen points. Such a sequence can be constructed using
a countable base for the topology of the coordinate polydisc,
choosing a new point at each step. After any finite number
of points has been removed, the sequence remains dense in $K$.

Every tail is Zariski dense in $W$. Indeed, if a polynomial
$q$ vanishes on a tail, then by continuity it vanishes on
$K$, and hence on a nonempty open part of the smooth
irreducible set $W$. Therefore, $q$ vanishes identically
on $W$.

\emph{3. A compact cube and normal convergence.}
Set $\delta_r=2^{-r}$ and
\[
 X=\prod_{r=1}^\infty[\delta_r/2,\delta_r].
\]
This is a nonempty compact complete metric space with metric
\[
 \rho(c,\widetilde c)=\sum_{r=1}^\infty
 2^{-r}\frac{|c_r-\widetilde c_r|}{\delta_r}.
\]
For $c\in X$, define
\begin{equation}\label{pm:polarized-series}
 U_c(z)=\sum_{r=1}^\infty c_r b(\lambda_r)e^{\lambda_r\cdot z}.
\end{equation}
On every compact set, uniformly in $c\in X$, this series
and all its derivatives converge normally. For a multi-index
$\alpha$, we have
\[
 \left\|c_r b(\lambda_r)\lambda_r^\alpha
 e^{\lambda_r\cdot z}\right\|
 \le\delta_rBT^{|\alpha|}e^{T\|z\|},
 \qquad\sum_r\delta_r=1.
\]
Thus $U_c$ is entire and
\begin{equation}\label{pm:type-bound}
 \|U_c(z)\|\le Be^{T\|z\|}.
\end{equation}
Each term of the series satisfies the system
by~(\ref{pm:polarization}); normal convergence of the derivatives
gives $A(D)U_c=0$. The normalization $b_j=1$ yields
\begin{equation}\label{pm:scalar-component}
 u_c(z):=(U_c)_j(z)=\sum_{r=1}^\infty c_r e^{\lambda_r\cdot z}.
\end{equation}

\emph{4. Independent finite coordinate systems on every tail.}
Let $\mathcal T_W$ denote the space of all finite sums of
exponentials with frequencies on $W$. For an integer
$q\ge0$, consider the space of functions on $W$
\[
 V_q=\Span_\C\{\lambda^\alpha:|\alpha|\le q\}\subset\C[W],
 \qquad h_q=\dim_\C V_q.
\]
Choose a basis $f_1,\ldots,f_{h_q}$ of this space.
For every $L$, the vectors
\[
 (f_1(\lambda_r),\ldots,f_{h_q}(\lambda_r)),\qquad r>L,
\]
span $\C^{h_q}$. Otherwise, a nonzero linear combination
of the $f_i$ would vanish on the entire tail, and hence
identically on $W$, contradicting the choice of a basis.
Consequently, there are indices $r_1,\ldots,r_{h_q}>L$
for which the matrix $(f_i(\lambda_{r_a}))$ is nonsingular.

This means that variations of the coefficients at these finitely
many frequencies fill the entire space of admissible $q$-jets
of test solutions at the origin. At another fixed point,
each column of the matrix is multiplied by a nonzero exponential,
and the conclusion remains valid. In particular, restricting
the frequencies to the compact set $K$ does not reduce
any of the finite-dimensional spaces of test jets.

\emph{5. Excluding a single differential relation.}
Let $\mathcal H_W$ be the set of all autonomous differential
polynomials with integer coefficients that are nontrivial on
$\mathcal T_W$. It is countable, since each polynomial
depends on finitely many derivatives and has finitely many
integer coefficients. For $F\in\mathcal H_W$, set
\[
 Z_F=\{c\in X:F[u_c](0)=0\}.
\]
The map $c\mapsto F[u_c](0)$ is continuous: it depends
on finitely many derivatives, and their series converge
uniformly in $c\in X$. Thus $Z_F$ is closed.

We show that $Z_F$ has empty interior. Nontriviality of
$F$ on $\mathcal T_W$ can be checked at the origin.
If $F[v](z_0)\ne0$ for some test solution, then
$v(z+z_0)$ also belongs to $\mathcal T_W$, and autonomy
gives a nonzero value at the origin. Thus, for sufficiently
large $q$, the polynomial $F$ defines a nonzero polynomial
on the space of admissible $q$-jets.

Take a nonempty open set $G\subset X$. It contains
a nonempty cylinder restricting only the first $L$
coordinates for some $L$. Fix admissible values of
these coordinates and all the remaining coefficients
except those chosen in step~4,
$c_{r_1},\ldots,c_{r_{h_q}}$, whose indices are greater
than $L$. All such variations remain in the cylinder.

The jet of the fixed part of the infinite series belongs
to the same space of admissible jets: it is a limit of
jets of finite test sums, and this space is finite-dimensional
and closed. The chosen free coefficients define a linear
isomorphism onto the whole space of independent jet
coordinates. Hence $F[u_c](0)$ is a nonzero polynomial
in these free complex coefficients. Adding a fixed
admissible jet does not change its nontriviality.

A nonzero complex polynomial cannot vanish on a product
of real intervals
\[
 \prod_{a=1}^{h_q}[\delta_{r_a}/2,\delta_{r_a}].
\]
This follows by successive application of the identity
theorem in each variable. Therefore, the chosen cylinder
contains a point $c$ with $F[u_c](0)\ne0$.
Thus $Z_F$ is closed and nowhere dense.

\emph{6. Applying Baire's theorem and excluding all finite
representations.} By Baire's theorem, the set
\begin{equation}\label{pm:coefficient-residual}
 X^*=X\setminus\bigcup_{F\in\mathcal H_W}Z_F
\end{equation}
is dense and is of type $G_\delta$; in particular, it is
nonempty. Take $c\in X^*$ and suppose that $u_c$ has a
finite representation $n\to k$ on some nonempty open set.

Apply Proposition~\ref{sc:elimination-proposition} to $W$.
The condition used here is precisely $\dim\Aff(W)>k$.
The affine hull of $W$ is not assumed to be all of
$\C^n$: in the lemma, the tangent field is chosen
outside the space
\[
 K_W=\{a\in\C^n:a\cdot\lambda\mbox{ is constant on }W\}.
\]
By Proposition~\ref{sc:elimination-proposition}, there is
an autonomous differential polynomial $F$ with integer
coefficients such that $F[u_c]=0$ on the domain of the
representation and which is nontrivial on the translated
test family $u_c+\mathcal T_W$.

In the present case, this also implies nontriviality on
$\mathcal T_W$. Indeed, $u_c$ belongs to the
closure of $\mathcal T_W$ in the topology of uniform convergence on compact sets,
by~(\ref{pm:scalar-component}).
At each fixed point, every finite jet of this function
belongs to the finite-dimensional linear space of test
jets. Translating this space by the jet of $u_c$
therefore leaves it unchanged. A polynomial in finitely
many derivatives is nontrivial on the translated family
if and only if it is nontrivial on the original test
family. Hence $F\in\mathcal H_W$.

The function $F[u_c]$ is entire. Its vanishing on a
nonempty open set extends to all of $\C^n$, whence
$F[u_c](0)=0$. This contradicts $c\in X^*$.
Thus $N_{n,k}(u_c)=\infty$.
Together with~(\ref{pm:type-bound}), this completes the proof.
\endproof

In Proposition~\ref{pm:exp-type}, residuality is established
in the compact coefficient cube $X$; all corresponding
solutions satisfy a common bound of exponential type.
In Theorem~\ref{pm:matrix-criterion}, residuality is established
in the whole Fr\'echet space $\mathcal E(A)$.
These are two separate assertions. They do not imply
that the subset of solutions of exponential type is
residual in the whole space $\mathcal E(A)$.

\subsection{An explicit vector family in three variables}
\label{ex:vector}

Consider the system for two unknown functions
\[
 v=u_x,\qquad u_y=v_x,\qquad u_z=v_{xx}.
\]
Its matrix of full symbols is
\[
 A(\xi,\eta,\zeta)=
 \left(\begin{array}{cc}
 -\xi&1\\
 \eta&-\xi\\
 \zeta&-\xi^2
 \end{array}\right).
\]
The condition $\rank A<2$ is equivalent to
$\eta=\xi^2$, $\zeta=\xi^3$. Thus,
\[
 W=\{(t,t^2,t^3):t\in\C\},\qquad
 b(t)=\left(\begin{array}{c}1\\t\end{array}\right)
\]
define the characteristic curve and a normalized polarization.
The affine hull of the curve is three-dimensional.

Choose the frequencies at $t_r=1/r$ and set
\begin{equation}\label{ex:vector-series}
 \left(\begin{array}{c}u_c(x,y,z)\\v_c(x,y,z)\end{array}\right)
 =\sum_{r=1}^\infty c_r
 \left(\begin{array}{c}1\\r^{-1}\end{array}\right)
 \exp\left(\frac xr+\frac y{r^2}+\frac z{r^3}\right),
 \qquad 2^{-r-1}\le c_r\le2^{-r}.
\end{equation}
Every tail of frequencies is Zariski dense in $W$:
the restriction of a polynomial to this curve is a
polynomial in $t$, which cannot vanish at all $t_r$
in a tail unless it is identically zero. Thus the proof
of Proposition~\ref{pm:exp-type} also applies to this
sequence; density of the frequencies in the whole
coordinate polydisc was only a convenient sufficient
condition for Zariski density of every tail.

All series~(\ref{ex:vector-series}) satisfy the same bound
\[
 \|(u_c,v_c)\|\le\sqrt2\,e^{|x|+|y|+|z|}
 \le\sqrt2\,e^{\sqrt3\,\|(x,y,z)\|}.
\]
For a dense $G_\delta$ set of coefficient sequences,
the function $u_c$ admits no finite representation
$3\to2$ on any nonempty open domain. The assertion
selects the coefficients by a Baire category argument;
it does not verify this property for a particular
arbitrarily prescribed sequence $c_r$ in the cube.

\subsection{Harmonic solutions of exponential type}
\label{ex:harmonic}

The parametrization
\[
 \lambda(t)=(t-t^3,i(t+t^3),2t^2)
\]
lies on the characteristic cone of the Laplace equation:
\[
 (t-t^3)^2+(i(t+t^3))^2+(2t^2)^2=0.
\]
Its image is a closed algebraic curve: the parameter
$t=(\lambda_1-i\lambda_2)/2$ is recovered linearly
from the image coordinates. The functions
$1,t-t^3,i(t+t^3),2t^2$ are linearly independent,
so the affine hull of this curve is $\C^3$.

For $t_r=1/r$ and $2^{-r-1}\le c_r\le2^{-r}$,
the series
\[
 h_c(x,y,z)=\sum_{r=1}^\infty
 c_r\exp((t_r-t_r^3)x+i(t_r+t_r^3)y+2t_r^2z)
\]
are entire solutions of $h_{xx}+h_{yy}+h_{zz}=0$
satisfying a common bound of exponential type.
Every tail of the chosen frequencies is Zariski dense
in the indicated curve. The coefficient cube argument
therefore yields a dense $G_\delta$ set of sequences
for which $h_c$ has infinite complexity $3\to2$
on every nonempty open subset. Here it suffices to
take, in the scalar construction, an irreducible
subset of the spectrum whose affine hull has
dimension greater than $k$; this subset need not
be a whole characteristic component.

\section{Topology of superposition classes}\label{top:section}

In this section we distinguish between two assertions: that the set of solutions
of infinite complexity contains a dense $G_\delta$ set, and that this set itself is a dense
$G_\delta$ set. The first already follows from autonomous polynomial obstructions.
For the second, we use normal families of internal functions. All spaces of entire
functions are endowed with the topology of uniform convergence on compact sets.

Operations on superposition representations require a common domain of definition.
If $A$ and $B$ have representations of depths $a$ and $b$ on the same open set,
then on this set
\[
N_{n,k}(A+B)\le 1+\max\{a,b\}.
\]
The same bound holds for the difference. For $k\ge2$, the product is constructed
similarly; for $k=1$, one can use the identity $4AB=(A+B)^2-(A-B)^2$. However,
if $N_{n,k}$ is defined by taking the minimum over arbitrary nonempty open subsets,
a common domain for two minimal representations cannot be assumed without justification.
In the proofs of the spectral assertions, such a domain exists because the exponential
being added has the required representation near every point. Below, a common domain
is fixed as part of the definition of the representation classes themselves.

\subsection{A countable normalization of representations}

Let $\mathcal F_{\le\ell}$ denote the set of entire functions on $\C^n$ admitting
a representation of depth at most $\ell$ on some nonempty open set. Put
\[
\mathcal F=\bigcup_{\ell\ge0}\mathcal F_{\le\ell}.
\]
Thus, $\mathcal F$ consists of functions of finite complexity on at least one open
subset. A combinatorial scheme specifies the arrangement of vertices, their arities,
and the choice of original coordinates. For $k=1$, the additions allowed by the
convention can be regarded as separate fixed vertices; normal families will be needed
only for the arbitrary analytic functions.

\begin{proposition}\label{top:fsigma}
The sets $\mathcal F_{\le\ell}$ and $\mathcal F$ are $F_\sigma$ sets in
$\mathcal O(\C^n)$: each is a countable union of closed subsets.
Each closed class in this union guarantees a representation on its own fixed
nonempty polydisc.
\end{proposition}

\proof
Fix a finite scheme. At each vertex representing an arbitrary analytic function $g_v$
of $p_v$ arguments, choose an open polydisc
\[
D_v=\{z\in\C^{p_v}:|z_j-c_{v,j}|<R_{v,j},\ 1\le j\le p_v\}
\]
and a closed inner polydisc
\[
K_v=\{z\in\C^{p_v}:|z_j-c_{v,j}|\le r_{v,j},\ 1\le j\le p_v\},
\qquad 0<r_{v,j}<R_{v,j}.
\]
All radii are chosen to be rational, and the centers have rational real and imaginary
parts. In particular, $K_v$ is compact and lies strictly inside $D_v$.
Also choose a positive integer $M_v$ and a nonempty polydisc $B$ in the original
variables, with rational parameters.

Let $\mathcal D$ denote this entire finite collection of data. Define
$\mathcal C_{\mathcal D}$ to be the class of entire functions $u$ for which there exist
functions $g_v$, holomorphic on $D_v$, such that
\[
|g_v(z)|\le M_v\quad(z\in D_v),\qquad s_v(X)\in K_v\quad(X\in B),
\]
where $s_v$ is the tuple of arguments of the corresponding vertex, and the output
of the scheme agrees with $u$ on $B$. All constraints here are non-strict.
No additional internal functions are introduced for the fixed addition vertices.

We show that $\mathcal C_{\mathcal D}$ is closed. Let $u_\nu\in\mathcal C_{\mathcal D}$
and $u_\nu\to u$ uniformly on compact subsets of $\C^n$. For each $u_\nu$, choose
one admissible family of internal functions $g_{v,\nu}$. By Montel's theorem, from
the bounded family on each $D_v$ one can extract a subsequence converging uniformly
on compact sets to a holomorphic function $g_v$. There are finitely many vertices,
so a single subsequence gives convergence simultaneously at all vertices.
The bounds $|g_v|\le M_v$ are preserved.

Proceed from the leaves to the root. Assuming convergence of the values at the
descendants, the arguments $s_{v,\nu}$ converge uniformly on compact subsets of $B$.
Their limits belong to $K_v$, since $K_v$ is closed. The convergence
$g_{v,\nu}\to g_v$ is uniform on $K_v$ and on some compact neighborhood of it
inside $D_v$. It follows that
\[
g_{v,\nu}(s_{v,\nu}(X))\longrightarrow g_v(s_v(X))
\]
uniformly on compact subsets of $B$. The same convergence is immediate for the fixed
addition vertices. Induction yields a limiting representation whose output equals
$u$ on $B$. Hence $u\in\mathcal C_{\mathcal D}$. The space of entire functions is
metrizable, so the sequential closedness just established is equivalent to closedness.

We now verify that these classes cover all local representations with the given
scheme. Choose a point in the domain of a representation and a sufficiently small
neighborhood of it. For each vertex, one can choose a polydisc $D_v$ of the specified
rational form so that its closure lies in the domain of holomorphy of the corresponding
internal function, while the value of the arguments at the chosen point lies strictly
inside some $K_v$. After shrinking the original neighborhood, all argument values
lie in these $K_v$. Choose a nonempty polydisc $B$ with rational parameters inside
this neighborhood. The internal functions are bounded on the closures of $D_v$,
so suitable integers $M_v$ exist. Thus the representation belongs to one of the
classes $\mathcal C_{\mathcal D}$.

There are countably many collections of rational data and integer bounds.
At a fixed depth there are finitely many combinatorial schemes, and there are
countably many finite schemes in total. Each class $\mathcal C_{\mathcal D}$ does
consist of locally representable functions, and the covering just verified gives
equality with the union of these closed classes. This proves the assertion for both
$\mathcal F_{\le\ell}$ and $\mathcal F$.
\endproof

The use of closed polydiscs $K_v$ is essential for the actual closedness of an
individual fixed class. With strict inequalities alone, a limit could reach the
boundary of the allowed set of arguments. The conditions
$|s_{v,j}-c_{v,j}|\le r_{v,j}<R_{v,j}$ are both preserved under passage to the limit
and leave a positive distance to the boundary of the domain of holomorphy $D_v$.

\subsection{The full set of solutions of infinite complexity}

\begin{theorem}\label{top:whole-gdelta}
Let $E$ be a closed complex linear subspace of $\mathcal O(\C^n)$ containing
a function $v$ that has no finite representation of the type under consideration
on any nonempty open set. Then the set
\[
E\setminus\mathcal F
\]
itself is a dense $G_\delta$ set in $E$.
\end{theorem}

\proof
By Proposition~\ref{top:fsigma}, we have
$E\cap\mathcal F=\bigcup_{\mathcal D}(E\cap\mathcal C_{\mathcal D})$, where the union
is countable and all its members are closed in $E$. We prove that each has empty
interior.

Suppose that some nonempty open set $U\subset E$ is contained in
$\mathcal C_{\mathcal D}$. Choose $u\in U$. By continuity of the map
$t\mapsto u+tv$, the functions $u+tv$ belong to $U$ for all sufficiently small
$t\in\C$. Fix such a $t\ne0$. The functions $u$ and $u+tv$ have representations
on the same polydisc $B$, fixed in the data $\mathcal D$. Their difference also
has a finite representation on $B$, since addition, multiplication by a constant,
and finite composition are allowed. Thus
\[
v=t^{-1}((u+tv)-u)
\]
has a finite representation on $B$, contrary to the choice of $v$.
Consequently, each $E\cap\mathcal C_{\mathcal D}$ is nowhere dense.

The closed subspace $E$ of the Fr\'echet space of entire functions is itself
a Fr\'echet space and, in particular, a Baire space. The complement of a countable
union of closed nowhere dense sets is dense and is a $G_\delta$ set.
This complement is precisely $E\setminus\mathcal F$.
\endproof

\begin{corollary}\label{top:system-gdelta}
If some irreducible component $W\subset V(I)$ satisfies
$\dim\Aff(W)>k$, then the full set of entire solutions of the system having no finite
representation of type $n\to k$ on any nonempty open subset is a dense
$G_\delta$ set in $\mathcal E(I)$.
\end{corollary}

\proof
Theorem~\ref{sc:criterion} gives an entire solution $v$ with the required property.
The space $\mathcal E(I)$ is a closed linear subspace of $\mathcal O(\C^n)$, since
it is the common kernel of continuous linear differential operators with constant
coefficients. Apply Theorem~\ref{top:whole-gdelta}.
\endproof

\begin{corollary}\label{top:matrix-gdelta}
For the index $j$ chosen in Theorem~\ref{pm:matrix-criterion}, the full set
\[
\{U\in\mathcal E(A):N_{n,k}(u_j)=\infty\}
\]
is a dense $G_\delta$ set in $\mathcal E(A)$.
\end{corollary}

\proof
By Corollary~\ref{top:system-gdelta}, the corresponding scalar set is a dense
$G_\delta$ set in $\mathcal E(I_j)$. Apply the continuous open surjection
$\pi_j:\mathcal E(A)\to\mathcal E(I_j)$ from Lemma~\ref{pm:projection}.
Its inverse image is precisely the displayed set.
\endproof

This topological argument strengthens the genericity assertion, but is not needed
for the existence of solutions of infinite complexity or for the qualitative criterion.
Their proof by means of autonomous polynomial obstructions already gives a dense
$G_\delta$ set consisting of such solutions. The normalization using Montel's theorem
is needed only to establish that the full set of solutions admitting no finite
representation anywhere is a $G_\delta$ set.

\section{The depth spectrum and quantitative bounds}\label{sp:section}

In this section we retain the numbering of classes adopted above. For $k\ge2$,
the zeroth class consists of functions of at most $k$ original coordinates, and each
subsequent level adds one composition with a function of at most $k$ arguments.
For $k=1$, the first level has the form $g(a_1(x_1)+\cdots+a_n(x_n))$, and subsequent
levels are constructed using binary operations. The numerical values of depth depend
on this convention. We define the spectrum of entire solutions by
\[
\mathcal S_{n,k}(I)=\{N_{n,k}(u):u\in\mathcal E(I)\}.
\]
Here, as above, complexity means the least depth of a representation on some nonempty
open set; the value $\infty$ means that no such representation exists on any nonempty
open set.

\subsection{Unbounded depths of finite sums of exponentials}

For an irreducible algebraic set $W\subset V(I)$, let $\mathcal T_W$ denote the space
of finite sums
\[
v(X)=\sum_{\nu=1}^q c_\nu\exp(p_\nu\cdot X),\qquad
p_\nu\in W,\quad c_\nu\in\C.
\]
All these functions are entire solutions of the original system.

\begin{lemma}\label{sp:unbounded-exponentials}
If $\dim\Aff(W)>k$, then for every integer $\ell\ge0$ there exists $v\in\mathcal T_W$ with
\[
\ell<N_{n,k}(v)<\infty.
\]
\end{lemma}

\proof
By Corollary~\ref{sc:level-polynomial}, there exists an autonomous polynomial
$F_\ell$ that gives a necessary condition for elements of $\mathcal T_W$ of depth
at most $\ell$ and is nontrivial on $\mathcal T_W$. Its construction uses the
finiteness of the list of output polynomials at a fixed depth, established in
Lemma~\ref{sc:finite-list}. Mere countability of the polynomials would not suffice
for this strengthening.

There is a finite sum $v\in\mathcal T_W$ for which $F_\ell[v]$ is not identically zero.
If necessary, translate the independent variables to obtain $F_\ell[v](0)\ne0$;
translation preserves $\mathcal T_W$. If $v$ had depth at most $\ell$ on some open set,
then $F_\ell[v]=0$ would hold there. By the identity theorem, this identity would
hold throughout $\C^n$, contrary to the choice of $v$. The finiteness of the complexity
of $v$ follows from its explicit construction as a finite sum of exponentials.
\endproof

\subsection{The exact spectrum and its initial threshold}

For $p\in\C^n$, let $a(p)$ denote the number of nonzero coordinates of $p$.
For $k\ge2$, a tree of depth $\ell$ has at most $k^\ell$ leaves, each depending
on at most $k$ original coordinates. Hence every function of this depth depends
on at most $k^{\ell+1}$ original coordinates.

\begin{lemma}\label{sp:exponential-depth}
Let $k\ge2$. If $a(p)\le k$, then $N_{n,k}(\exp(p\cdot X))=0$. If $a(p)>k$, then
\[
N_{n,k}(\exp(p\cdot X))=\lceil\log_k a(p)\rceil-1.
\]
Under the convention adopted for $k=1$, every exponential has complexity at most one,
and for $a(p)\ge2$ its complexity equals one.
\end{lemma}

\proof
The exponential depends essentially on each coordinate for which $p_i\ne0$.
Counting leaves gives the lower bound for $k\ge2$. For the upper bound, partition
the nonzero coordinates into blocks of size at most $k$. The exponential of the sum
in each block belongs to the zeroth class. Multiplying these exponentials in
a balanced $k$-ary tree requires the least depth $\ell$ for which $k^{\ell+1}\ge a(p)$.
For $k=1$, the representation $\exp(\sum_i p_ix_i)$ is at the first level, whereas
a function of the zeroth class cannot depend essentially on two original coordinates.
\endproof

\begin{theorem}\label{sp:spectrum}
Let $W\subset V(I)$ be irreducible and satisfy $\dim\Aff(W)>k$.

For $k=1$, we have the exact equality
\[
\mathcal S_{n,1}(I)=\{0,1,2,\ldots,\infty\}.
\]
For $k\ge2$, put
\[
a_W=\#\{i:\xi_i\mbox{ is not identically zero on }W\},\qquad
e_W=\lceil\log_k a_W\rceil-1.
\]
Then $a_W>k$ and
\[
\{0\}\cup\{e_W,e_W+1,e_W+2,\ldots\}\cup\{\infty\}\subset\mathcal S_{n,k}(I).
\]
If $a_W\le k^2$, in particular if $n\le k^2$, the spectrum is full:
\[
\mathcal S_{n,k}(I)=\{0,1,2,\ldots,\infty\}.
\]
In all the assertions above, each realizable finite value in the stated spectrum
is attained by a finite sum of characteristic exponentials.
\end{theorem}

\proof
The inequality $a_W\ge\dim\Aff(W)>k$ is immediate. At a generic point, all
$a_W$ coordinates that are not identically zero on $W$ are nonzero: a finite union
of proper closed subsets does not cover an irreducible set $W$.
By Lemma~\ref{sp:exponential-depth}, the corresponding exponential has depth exactly
$e_W$, and all exponentials with frequencies in $W$ have depth at most $e_W$.

Fix an integer $L>e_W$. By Lemma~\ref{sp:unbounded-exponentials}, there is a finite sum
\[
v_q=\sum_{j=1}^q c_j\exp(p_j\cdot X),\qquad N_{n,k}(v_q)>L.
\]
Put $v_0=0$ and let $v_j$ denote the corresponding partial sums. On adding one term,
we have
\[
N_{n,k}(v_{j+1})\le1+\max\{N_{n,k}(v_j),e_W\}.
\]
There is no difficulty concerning a common domain for the representations: the
exponential has the stated representation near every point, in particular on the
domain where a representation of $v_j$ has been chosen. Choose the first index $j$
for which $N_{n,k}(v_j)\ge L$. Then $N_{n,k}(v_{j-1})\le L-1$ and $e_W\le L-1$,
whence $N_{n,k}(v_j)\le L$. Therefore, $N_{n,k}(v_j)=L$. This argument does not
require the complexities of the partial sums to be monotone.

The value $e_W$ is realized by a single exponential, and the value zero by the zero
function. The value $\infty$ is supplied by Theorem~\ref{sc:criterion}.
If $a_W\le k^2$, then $e_W=1$ and the resulting spectrum is full.

For $k=1$, the same argument applies with $e_W=1$. There is a point of $W$ with
at least two nonzero coordinates, since $\dim\Aff(W)>1$; the corresponding exponential
has complexity exactly one.
\endproof

\begin{corollary}\label{sp:three-two-spectrum}
In the case $3\to2$, failure of the uniform finiteness condition in
Theorem~\ref{sc:criterion} always gives the full spectrum $\{0,1,2,\ldots,\infty\}$.
The same holds for $n\to n-1$ when $n\ge3$.
\end{corollary}

\proof
For $3\to2$, we have $3\le2^2$. For $n\ge3$, we have $n\le(n-1)^2$.
Apply Theorem~\ref{sp:spectrum}.
\endproof

\subsection{A system with an initial gap in the spectrum}

\begin{proposition}\label{sp:gap}
Let $k\ge2$ and $n>k^2$. Consider the curve
\[
W=\{(t,t^{-1},t^2,t^3,\ldots,t^{n-1}):t\in\C,\ t\ne0\}
\]
and its ideal
\[
J=(\xi_1\xi_2-1,\ \xi_j-\xi_1^{j-1}\ (3\le j\le n)).
\]
If
\[
e=\lceil\log_k n\rceil-1,
\]
then the exact spectrum of entire solutions is
\[
\mathcal S_{n,k}(J)=\{0,e,e+1,e+2,\ldots,\infty\}.
\]
In particular, the values $1,\ldots,e-1$ are not realized.
\end{proposition}

\proof
The quotient ring $\C[\xi_1,\ldots,\xi_n]/J$ is isomorphic to $\C[t,t^{-1}]$.
Thus $W$ is a closed irreducible affine curve with ideal $J$. The functions
$1,t,t^{-1},t^2,\ldots,t^{n-1}$ are linearly independent, so $\Aff(W)=\C^n$.

In addition to the first displayed generator, $J$ contains $\xi_j\xi_2^{j-1}-1$
for $3\le j\le n$. Consequently, every solution satisfies
\[
D_1D_2u=u,\qquad D_jD_2^{j-1}u=u\quad(3\le j\le n).
\]
If a solution is independent of $x_1$ or $x_2$, the first equation gives $u=0$.
If it is independent of $x_j$, $j\ge3$, the corresponding second equation also
gives $u=0$. Thus every nonzero solution depends essentially on all $n$ original
coordinates on every nonempty open subset.

Counting leaves rules out depths smaller than $e$ for nonzero solutions. All
frequency coordinates on $W$ are nonzero, so $e_W=e$. Theorem~\ref{sp:spectrum}
realizes every value starting at $e$, as well as $\infty$. The zero function has
complexity zero.
\endproof

\subsection{Thickening a coordinate $k$-plane}

Put $c=n-k$ and, for an integer $m\ge1$, consider the ideal
\[
I_m=(\xi_{k+1},\ldots,\xi_n)^m.
\]
Use coordinates $X=(s,t)$, where $s\in\C^k$ and $t\in\C^c$.
Let $h_{n,k}(m)$ denote the supremum of the complexities of all local solutions
of this system, and put
\[
M_m={m+c-1\choose c}.
\]

\begin{theorem}\label{sp:quantitative}
For $k\ge2$ and $m\ge1$, the following bounds hold:
\[
\lceil\log_k((k-1)M_m+1)\rceil-1\le h_{n,k}(m)
\le\lceil\log_k n\rceil-1+\lceil\log_k M_m\rceil.
\]
Moreover, $h_{n,k}(1)=0$, and for fixed $n,k$ as $m\to\infty$,
\[
h_{n,k}(m)=(n-k)\log_k m+O_{n,k}(1).
\]
For $k=1$, with the adopted convention for the first level,
\[
h_{n,1}(m)=(n-1)\log_2 m+O_n(1),\qquad h_{n,1}(1)=0.
\]
\end{theorem}

\proof
The equations $D_t^\beta u=0$, $|\beta|=m$, show that the general local solution has the form
\[
u(s,t)=\sum_{|\beta|<m}t^\beta f_\beta(s),
\]
where $f_\beta$ are arbitrary holomorphic functions of $k$ variables.
The number of these functions is $M_m$.

First let $k\ge2$. To construct one summand, take one leaf $f_\beta(s)$ and partition
the transverse coordinates into blocks of size at most $k$. The monomial associated
with each block is also an admissible leaf. In total, at most
\[
1+\left\lceil\frac{n-k}{k}\right\rceil=\left\lceil\frac nk\right\rceil
\]
leaves are needed. Their product is constructed by a balanced $k$-ary tree of depth
at most $\lceil\log_k n\rceil-1$. Summing the $M_m$ terms adds at most
$\lceil\log_k M_m\rceil$ levels. This gives the upper bound. For $m=1$, the solution
has the form $f(s)$ and already belongs to the zeroth class.

For the lower bound, consider the space of realizable $r$-jets of solutions at the
origin. For $r\ge m-1$, its dimension is
\[
d_m(r)=\sum_{|\beta|<m}{r-|\beta|+k\choose k}
=\frac{M_m}{k!}r^k+O(r^{k-1}).
\]
Each such jet is realized by a polynomial solution: it suffices to choose the
corresponding finite Taylor expansions of the functions $f_\beta$.

A tree of depth $\ell$ has at most
\[
T_\ell=\frac{k^{\ell+1}-1}{k-1}
\]
vertices, including the leaves. Its output $r$-jet is determined by at most
$T_\ell{r+k\choose k}$ Taylor coefficients of the internal functions.
At each vertex, the series is centered at the value of its internal arguments.
The truncated composition rule is polynomial in these coefficients. At a fixed
depth, there are finitely many schemes.

If $M_m>T_\ell$, then for all sufficiently large $r$ the dimension of the space
of solution jets exceeds the dimension of the Zariski closure of the jets of
functions of depth at most $\ell$. Hence there is a polynomial $F$ in finitely many
derivatives that vanishes on all jets of this class but is not identically zero on
the space of solution jets. The relation is autonomous: translation of the original
point is absorbed by recentering the arbitrary internal functions.

Choose a polynomial solution $u$ with $F[u](0)\ne0$. A representation of depth
at most $\ell$ on any other open set would give $F[u]=0$ there, and then everywhere
by the identity theorem. Consequently, $N_{n,k}(u)>\ell$. Optimizing the condition
$M_m>T_\ell$ over integers $\ell\ge0$ gives
\[
h_{n,k}(m)\ge\lceil\log_k((k-1)M_m+1)\rceil-1.
\]
The relation $M_m=m^c/c!+O(m^{c-1})$ gives both asymptotic bounds for $k\ge2$.

For $k=1$, the number of internal functions of one variable at level $\ell\ge1$
does not exceed
\[
A_\ell=(n+2)2^{\ell-1}-1.
\]
Indeed, the first level contains $n+1$ such functions, and subsequent levels satisfy
the recurrence $A_{\ell+1}=2A_\ell+1$. The dimension of the parameter space for $r$-jets
does not exceed $A_\ell(r+1)$. Comparison with $d_m(r)=M_mr+O(1)$ when $M_m>A_\ell$
gives the lower bound $\log_2 M_m-O_n(1)$.

For the upper bound, each nonzero term $t^\beta f_\beta(s)$, where $s$ is now a single
coordinate, has complexity at most one on a common regular domain:
\[
t^\beta f_\beta(s)=\exp\left(\log f_\beta(s)+\sum_{j=1}^{n-1}\beta_j\log t_j\right).
\]
The zero coefficient functions are discarded. For the finitely many remaining
functions, one can choose a common nonempty domain containing none of their zeros
and satisfying $t_j\ne0$, and then choose local branches of the logarithms.
This is allowed by the definition of complexity used here, on some nonempty open
set. Binary summation gives the upper bound $1+\lceil\log_2 M_m\rceil$.
Since $M_m={m+n-2\choose n-1}$, we obtain $h_{n,1}(m)=(n-1)\log_2 m+O_n(1)$.
\endproof

\begin{corollary}\label{sp:three-two-bounds}
For $I_m=(\xi_3^m)$ and complexity of type $3\to2$, we have
\[
\lceil\log_2(m+1)\rceil-1\le h_{3,2}(m)\le1+\lceil\log_2 m\rceil,
\qquad h_{3,2}(1)=0.
\]
The lower bound equals $\lfloor\log_2 m\rfloor$. These are two-sided bounds;
in general, they do not determine the exact values for small multiplicities.
\end{corollary}

\proof
In this case, $c=1$, $M_m=m$, and $\lceil\log_2 3\rceil-1=1$.
Apply Theorem~\ref{sp:quantitative}.
\endproof

In the case $3\to1$ and $I_m=(\xi_2,\xi_3)^m$, we obtain
$h_{3,1}(m)=2\log_2 m+O(1)$. Here
\[
u=\sum_{i+j<m}x_2^ix_3^jf_{ij}(x_1),\qquad
d_m(r)={m+1\choose2}(r+1)-\frac{m(m-1)(m+1)}3,\qquad r\ge(m-1).
\]
Thus multiplicities substantially change the maximal depth, although the radical
of the ideal remains unchanged.

\subsection{A single equation and examples in three variables}

\begin{corollary}\label{sp:one-equation}
Let $P\in\C[\xi_1,\ldots,\xi_n]$ be a nonzero nonconstant polynomial. If $k<n-1$,
the equation $P(D)u=0$ has entire solutions of infinite complexity of type $n\to k$,
and the depths of its finite characteristic exponential sums are unbounded.
If $k=n-1$, all local solutions have uniformly bounded complexity if and only if
$P$ factors into affine linear factors. The multiplicities of the factors are
arbitrary. If the latter condition fails, the spectrum is full.
\end{corollary}

\proof
Every irreducible component of $V(P)$ has dimension $n-1$. Thus, for $k<n-1$,
the dimension of its affine hull exceeds $k$. For $k=n-1$, such a component is
contained in an affine subspace of dimension at most $n-1$ if and only if it is
itself an affine hyperplane. For an irreducible factor of $P$, this is equivalent
to its linearity. Apply Theorem~\ref{sc:criterion}. Fullness of the spectrum in
the unbounded case follows, for $n=2$, from the case $k=1$, and, for $n\ge3$,
from Corollary~\ref{sp:three-two-spectrum}.
\endproof

In the following examples, the original variables are denoted by $x,y,z$ and the
frequency variables by $\xi,\eta,\zeta$.

\begin{enumerate}
\item The equation $u_z=0$ has characteristic plane $\zeta=0$.
Its solutions $f(x,y)$ have complexity zero of type $3\to2$.
For $3\to1$, the spectrum is full.

\item The equation $u_{zz}=0$ has the same characteristic plane, with multiplicity.
Its solutions have the form $f(x,y)+zg(x,y)$ and complexity at most two of type
$3\to2$. General multiplicity bounds are given in Corollary~\ref{sp:three-two-bounds}.

\item For the system $u_y=u_{xx}$, $u_z=0$, the characteristic set is the plane
parabola $(t,t^2,0)$. All solutions have complexity zero of type $3\to2$,
whereas the spectrum of complexities of type $3\to1$ is full.

\item For the system $u_y=u_{xx}$, $u_z=u_{xxx}$, the characteristic curve has the form
$(t,t^2,t^3)$. Its affine hull is $\C^3$, since $1,t,t^2,t^3$ are linearly independent.
Thus the spectrum is full for both $3\to1$ and $3\to2$.

\item For the Laplace equation $u_{xx}+u_{yy}+u_{zz}=0$, the characteristic set
is the irreducible quadric $\xi^2+\eta^2+\zeta^2=0$ with affine hull $\C^3$.
Irreducibility also follows from the fact that a quadratic form of rank three cannot
be a product of two linear forms. The spectrum is full for both $3\to2$ and $3\to1$;
in particular, there exist entire complex harmonic functions of infinite complexity
of both types.

\item The equation $u_z=u_{xx}+u_{yy}$ has characteristic surface
$\zeta=\xi^2+\eta^2$. It is irreducible and is not contained in an affine plane.
Hence its affine hull is $\C^3$, and both spectra, $3\to1$ and $3\to2$, are full.
\end{enumerate}

In all examples where the uniform finiteness condition fails, solutions
of infinite complexity can be chosen from a dense residual subset of the space of
entire solutions. The precise strengthening to an assertion about the full set of
solutions of infinite complexity is given in Section~\ref{top:section}.

\section{Zero-dimensional spectra and explicit formulas for thickened lines}
\label{lp:finite}

\subsection{A zero-dimensional scalar spectrum}

Suppose that $V(I)$ is finite and nonempty. Then the algebra $B=R/I$ is finite-dimensional.
Let $T_j$ denote multiplication by the class of $\xi_j$ in $B$. For $\lambda\in B^*$, set
\[
u_\lambda(x)=\lambda\left(\exp\left(\sum_{j=1}^n x_jT_j\right)[1]\right).
\]
This is an entire solution: differentiation with respect to $x_j$ corresponds to applying
$T_j$, and every element of $I$ acts as zero on $B$. Conversely, a local solution $u$
defines a functional $\lambda_u([P])=P(D)u(0)$. Equality of all Taylor coefficients
shows that these correspondences are mutually inverse. Hence
\[
\dim_{\C}S(I)=\dim_{\C}R/I,
\]
and every local solution extends to an entire one. Decomposing the finite-dimensional
space into joint generalized eigenspaces of the commuting operators $T_j$ gives
\[
u(x)=\sum_{p\in V(I)}e^{p\cdot x}q_p(x),
\]
where the degrees of the polynomials $q_p$ are bounded in terms of the ideal $I$.
If $V(I)=\emptyset$, then $I=R$ and the zero solution is the only solution.
For a finite-dimensional system module, similarly,
$\dim_{\C}S(M)=\dim_{\C}M$: a functional on $M$ recovers all components
by means of the same matrix exponentials.

\subsection{A thickened line in three variables}

After a permutation of coordinates, an affine line has equations
\[
\eta=a\xi+b,\qquad \zeta=c\xi+d_0.
\]
Let $\mathcal L=(\eta-a\xi-b,\zeta-c\xi-d_0)$. Under the substitution
\[
t=x+ay+cz,\quad s=y,\quad r=z,\qquad u=e^{bs+d_0r}v(t,s,r)
\]
the generators of $\mathcal L$ become $D_s,D_r$. Thus the general solution of
$\mathcal L^m(D)u=0$ is
\begin{equation}
u=e^{by+d_0z}\sum_{i+j<m}y^iz^jf_{ij}(x+ay+cz).
\label{lp:thick-line}
\end{equation}
The functions $f_{ij}$ are arbitrary; their number is ${m+1\choose 2}$. If
$\mathcal L^m\subset I$, the additional equations from $I$ impose relations on $f_{ij}$,
but preserve the upper bound on complexity. In two variables, the corresponding formula
contains $m$ functions of a single linear form.

For $k=1$ and the first level adopted here, a single exponential and a function of one
linear form have depth at most $1$. A coordinate monomial also has depth at most $1$
on a common regular domain: logarithms of the nonzero coordinates are used.
The depth bound concerns a representation on a nonempty open domain, so such a
restriction of the domain is allowed.

If the primary decomposition in three variables consists of components supported on
lines and points, choose exponents $m_\alpha,n_\gamma$ so that the corresponding powers
of the radicals are contained in the primary ideals. Set
\[
L=\sum_\alpha{m_\alpha+1\choose 2}+\sum_\gamma{n_\gamma+2\choose 3}.
\]
Every solution is a sum of at most $L$ terms, each of which is a product of at most
three functions of depth at most $1$. Successive application of the identity
$4AB=(A+B)^2-(A-B)^2$ and balanced summation give
\[
\sup_{u\in S(I)}N_{3,1}(u)\le 5+\lceil\log_2L\rceil\qquad (L>0).
\]
In two variables, the definition of $L$ uses $m_\alpha$ and
${n_\gamma+1\choose 2}$, respectively. This bound is not optimal, but explicitly
accounts for multiplicities.

\section{Analytic lifting of solutions from a projection}
\label{lp:projection}

Assume that $I\ne R$. We give here an additional description of local solutions of
systems in two and three variables when $\dim V(I)\le 1$. It is not needed to prove
the general criterion, but explains the conditions under which an arbitrary analytic
solution on a plane projection can be lifted.

\subsection{Choosing a noncharacteristic frequency coordinate}

Denote the frequencies by $(\xi,\eta,\zeta)$; in the two-dimensional case the third
coordinate is absent. Choose a linear coordinate $\xi$ that does not vanish at any
point at infinity of the projective closures of the one-dimensional components.
There are only finitely many such points, so this choice is possible. The independent
variables undergo the dual change, preserving the bilinear pairing with the
frequencies. Set $K=\C[\xi]$.

The ideal $I$ contains polynomials
\begin{equation}
Q_\eta=\eta^a+\sum_{i<a}a_i(\xi)\eta^i,\quad \deg a_i\le a-i;
\qquad
Q_\zeta=\zeta^b+\sum_{j<b}b_j(\xi)\zeta^j,\quad \deg b_j\le b-j.
\label{lp:monic}
\end{equation}
Indeed, the projection of the projective closure onto the plane $[w:\xi:\eta]$ is
defined on that closure: its center at infinity with $\xi=0$ does not meet the closure.
The image is closed and does not contain the point $[0:0:1]$. Consequently, the
defining polynomial of each one-dimensional plane component contains a nonzero
term $\eta^{\deg Q}$. For isolated points, add the factors $\eta-\eta_0$.
The product vanishes on $V(I)$; some power of it belongs to $I$ by the Nullstellensatz.
After normalization, we obtain the first relation in (\ref{lp:monic}). The second is
obtained in the same way. If the spectrum is zero-dimensional, products of linear
factors in $\eta$ and $\zeta$ give the same conclusion directly.

In particular, $R/I$ is finitely generated over $K$. The total-degree bounds in
(\ref{lp:monic}) are also essential: the orders of the operators $Q_\eta(D)$ and
$Q_\zeta(D)$ are $a$ and $b$, and the hyperplanes $y=0$ and $z=0$ are
noncharacteristic for the corresponding problems.

\subsection{A description in terms of one-dimensional analytic data}

Let $\mathcal O_1=\C\{x\}$, with $\xi f=f'$.

\begin{lemma}
\label{lp:hom}
With the above choice of coordinates, the map
\[
u\longmapsto\Phi_u,\qquad \Phi_u([P])=P(D)u|_{y=z=0}
\]
defines an isomorphism
\[
S(I)\simeq\Hom_K(R/I,\mathcal O_1).
\]
\end{lemma}

\proof
Injectivity follows from analyticity. For a given homomorphism $\Phi$, set
$f_{ij}(x)=\Phi([\eta^i\zeta^j])$, $0\le i<a$, $0\le j<b$.
First, for each $j$, use the Cauchy--Kovalevskaya theorem to solve
\[
Q_\eta(D_x,D_y)v_j=0,\qquad D_y^iv_j|_{y=0}=f_{ij}.
\]
Then solve
\[
Q_\zeta(D_x,D_z)u=0,\qquad D_z^ju|_{z=0}=v_j.
\]
In the second problem, $y$ is an additional analytic variable. Commutativity of the
operators and uniqueness for the second problem give $Q_\eta(D)u=0$.
The recurrences for $D_y^iD_z^ju|_{y=z=0}$ coincide with those for
$\Phi([\eta^i\zeta^j])$ in the module $R/I$. Hence
$P(D)u|_{y=z=0}=\Phi([P])$ for all $P$. If $P\in I$, then for every multi-index
$\alpha$ we have $D^\alpha P(D)u(0)=0$. Analyticity gives $P(D)u=0$.
In the two-dimensional case, a single Cauchy problem suffices.
\endproof

The module $\mathcal O_1$ is divisible: every equation $p(D_x)h=f$ with $p\ne 0$
has an analytic solution germ. Over the principal ideal domain $K=\C[\xi]$, this
implies injectivity. In the finitely generated case needed here, the extension can
be constructed explicitly. Given an inclusion $B\subset C$, adjoin to $B$ a single
generator $c$. The ideal $\{p\in K:pc\in B\}$ is either zero or generated by
$p_0\ne 0$. In the first case the image of $c$ is arbitrary. In the second, choose
$h$ with $p_0(D_x)h=\Phi(p_0c)$ and set $\Phi(b+pc)=\Phi(b)+p(D_x)h$.
The choice of $p_0$ ensures that this is well defined. Repetition gives an extension
to $C$.

Thus the diagonal embedding $R/I\hookrightarrow\bigoplus_\alpha R/Q_\alpha$
for a primary decomposition gives analytic additivity once again when the spectrum
has dimension at most one. A single frequency coordinate $\xi$ can be chosen
simultaneously for all components. The sum of the solution spaces is generally not
direct.

\subsection{Lifting from a plane projection and the limits of its applicability}

Let $C\subset\C^3$ be an irreducible curve, with coordinates chosen as above.
For the closure of its plane projection, with equation $Q(\xi,\eta)=0$, we have
\[
B=\C[\xi,\eta]/(Q)\hookrightarrow C_0=\C[\xi,\eta,\zeta]/\mathcal I_C.
\]
Both modules are finite over $K$. A solution germ of $Q(D)v=0$ defines a
$K$-homomorphism $B\to\mathcal O_1$. Extend it to $C_0$ and apply
Lemma~\ref{lp:hom}. The resulting solution satisfies $u(x,y,0)=v(x,y)$.

This assertion concerns the chosen noncharacteristic projection. Finiteness of an
arbitrary algebraic projection does not by itself ensure convergence of the formal
lift. For example, for the curve $(t,t^2,t^3)$, the system
\[
u_y=u_{xx},\qquad u_z=u_{xxx}
\]
formally recovers the solution from $f(x)$ by the formula
\[
u(x,y,z)=\sum_{j,l\ge 0}\frac{y^jz^l}{j!l!}f^{(2j+3l)}(x).
\]
This series may diverge even when $f$ is analytic. Moreover, the entire plane solution
\[
v(x,y)=\sum_{r\ge 1}e^{-r^{5/2}}e^{rx+r^2y}
\]
has no holomorphic lift near the origin on the original plane $z=0$.
The coefficient of $z^l$ in the putative lift at $x=y=0$ is
\[
\frac{1}{l!}\sum_{r\ge 1}e^{-r^{5/2}}r^{3l}.
\]
For $r$ of order $l^{2/5}$, the $l$th root of this expression is at least
$c\,l^{1/5}$ for some $c>0$. The radius of convergence is zero. If $t^3$ is chosen
as the distinguished frequency, the condition on the choice of coordinates holds;
for example, the projection $(t^3,t)$ has equation $\mu^3-\lambda=0$.

\section{Explicit elimination in the first class}
\label{lp:first}

To illustrate the procedure, consider $u=g(a(x)+b(y))$. Set
\[
H=a'u_y-b'u_x=0,\qquad A=u_xu_{yy}-u_yu_{xy},\qquad B=u_xu_y.
\]
Since $D_ya'=0$, we obtain
\[
u_yD_yH-u_{yy}H=b'A-b''B=0.
\]
Then $D_xb'=D_xb''=0$ gives
\begin{equation}
BD_xA-AD_xB=0.
\label{lp:first-relation}
\end{equation}
This is first derived at a generic point where the required factors are nonzero;
after expansion it is a polynomial identity, and hence extends to the whole domain.
If one of the original factors vanishes identically, the relation is checked
directly. When $u_xu_y\ne 0$, it is equivalent to
\[
D_xD_y\log(u_x/u_y)=0.
\]
This differential description of the first class agrees with~\cite{Bel2007}.
Here it is essential to check that the relation is nontrivial on the solution space
of another equation.

For a test sum $w$ with frequencies on the parabola $(t,t^2)$, set
$U_j:=Z_{t^j}[w]=\partial_x^jw$. Substituting $w$ for $u$, we obtain
$A[w]=U_1U_4-U_2U_3$, $B[w]=U_1U_2$. The left-hand side of
(\ref{lp:first-relation}), evaluated at $w$, is
\[
U_1^2U_2U_5-U_1U_2U_3^2-(U_1U_4-U_2U_3)(U_2^2+U_1U_3).
\]
The new coordinate $U_5$ occurs with the nonzero coefficient $U_1^2U_2$.
Independence of the finite frequency coordinates shows that this is a nontrivial
polynomial.

Multiplicities change the maximal depth even in small examples. For
$I_1=(\eta,\zeta)$, we have $u=f(x)$ and the maximum is $0$.
For $I_2=(\eta^2,\zeta)$, the solutions have the form $f(x)+yg(x)$, and the
maximum is $2$ under the $3\to 1$ convention. The upper bound follows from
a representation of $yg(x)$ in the first class on a regular domain and one addition.
The lower bound is furnished by $u=x+ye^x$, since
\[
D_xD_y\log(u_x/u_y)=\frac{e^{-x}}{(e^{-x}+y)^2}\ne 0.
\]
If this function of three variables belonged to the first class, fixing $z$ would
yield a representation in the first class in two variables, contradicting the computation.

\section{Examples of the transfer of complexity under changes of variables}

The complexity of solutions of nonlinear PDEs was studied, in particular,
in~\cite{Bel2021}, where solutions of the Burgers equation of depth at most one
were described. In the following example, we are interested in how a change in a
parameter of the equation affects the existence of solutions of infinite complexity.

\subsection{A jump in complexity in a one-parameter family}

Consider the family of nonlinear equations
\[
2vv_{xy}-v_xv_y=4\lambda v^2,\qquad \lambda\in\C.
\]
When $\lambda=0$, every solution has depth of type $2\to 1$ at most one at a
generic point. For every $\lambda\ne 0$, there are entire solutions of exponential
type having infinite complexity on every nonempty open subset.

Indeed, the substitution $v=u^2$ gives the identity
\[
2vv_{xy}-v_xv_y-4\lambda v^2=4u^3(u_{xy}-\lambda u).
\]
Near a point where $v\ne 0$, choosing a branch of the square root recovers the
linear equation. When $\lambda=0$, integration gives $u=f(x)+g(y)$, whence
$v=(f(x)+g(y))^2$. The identically zero solution belongs to the zeroth class.
When $\lambda\ne 0$, the spectrum of the linear equation is the irreducible
hyperbola $\xi\eta=\lambda$, whose affine hull is $\C^2$. By the results of
Section~\ref{pm:exp-section}, one can choose an entire solution $u$ of exponential
type and infinite complexity. Its square retains exponential type. Any finite
scheme for $u^2$ would, after restriction to an open subset where $u\ne 0$ and
application of a one-variable branch of the square root, give a finite scheme for $u$.

\subsection{Increase and decrease of depth under a single blow-up}

In the chart $\pi(s,t)=(s,st)$ of the blow-up of the plane, the function $y$ of
depth zero becomes $st$, of depth one. Here and below, depth is computed at a
generic point. The product of two nonconstant functions of different variables has
depth one by the formula $ab=\exp(\log a+\log b)$.

The reverse change in depth occurs for $P=x^2+xy+y^2$. The representation
$P=(x+y)^2-xy$ gives the upper bound two. For a function in the first class,
$P=G(a(x)+b(y))$, the ratio $P_x/P_y=a'(x)/b'(y)$ has separated variables.
However,
\[
\partial_x\partial_y\log(P_x/P_y)=\frac{6(x^2-y^2)}{(2x+y)^2(x+2y)^2}\not\equiv 0.
\]
Thus $P$ has depth two, whereas its total transform $P\circ\pi=s^2(1+t+t^2)$
has depth one.

\subsection{A harmonic solution near the exceptional divisor}

Take an entire harmonic solution $h$ from the family in Section~\ref{ex:harmonic},
choosing the coefficients so that its complexity of type $3\to 2$ is infinite on
every open subset. In a chart of the blow-up of the origin in $\C^3$, set
\[
V(s,p,q)=h(s,sp,sq),\qquad R=s\partial_s-p\partial_p-q\partial_q.
\]
For $s\ne 0$, differentiation with respect to the original coordinates corresponds
to the operators $s^{-1}R$, $s^{-1}\partial_p$, $s^{-1}\partial_q$.
Since $R(s^{-1})=-s^{-1}$, the Laplace equation takes the form
\[
L_\sigma V=0,\qquad L_\sigma=R^2-R+\partial_p^2+\partial_q^2.
\]
Both sides are holomorphic also at $s=0$. On the exceptional divisor,
$V(0,p,q)=h(0,0,0)$ is constant; nevertheless, $V$ has infinite complexity in every
neighborhood of every point of that divisor. Indeed, a finite representation on an
open subset could be restricted to $s\ne 0$. The inverse substitution
$(s,p,q)=(x,y/x,z/x)$ would then give a finite representation of $h$, since each
of the ratios depends on only two coordinates.

Combining the blow-up with squaring, we obtain an entire solution $w=V^2$ of
the nonlinear equation with polynomial coefficients
\[
2wL_\sigma w-(Rw)^2-w_p^2-w_q^2=0.
\]
This follows from the differentiation rule for a square. On any open subset where
$w$ has a finite representation, one could choose a smaller open subset with
$V\ne 0$ and recover $V$ by taking a square root. Hence $w$ also has infinite
complexity in every neighborhood of the divisor.

\section*{Concluding remarks}
\addcontentsline{toc}{section}{Concluding remarks}

The scalar and matrix criteria are expressed by a single geometric condition:
the affine hull of each irreducible component of the spectral support must have
dimension at most the number of arguments allowed for the internal functions.
The algebraic dimension of the component itself is insufficient for this purpose.
A cubic space curve already exhibits the difference.

The qualitative classification depends on the support. Multiplicities and the
module structure determine the solution space within that support and the
quantitative bounds. For a thickened coordinate plane, the leading term of the
maximal depth is logarithmic in the multiplicity. Exact values for small
multiplicities and optimal bounds for general primary modules remain separate
problems.

If the spectral condition fails, there is an index $j$ for which the set of all
entire solutions $U=(u_1,\ldots,u_m)$ with $N_{n,k}(u_j)=\infty$ is a dense
$G_\delta$ in the space $\mathcal E(A)$
(Corollary~\ref{top:matrix-gdelta}). No growth restrictions on solutions are imposed
in this assertion. In \S~\ref{pm:exp-section}, we also constructed a family of
solutions $U_c$ satisfying a common bound of exponential type. For every $c$ in a
dense $G_\delta$ subset of the coefficient space $X$, the component $(U_c)_j$
has infinite complexity.

The author used ChatGPT in the preparation of this article.


\begin{thebibliography}{99}
\small\raggedright
\setlength{\itemsep}{5pt}
\setlength{\parsep}{0pt}
\setlength{\parskip}{0pt}

\bibitem{Bel2007}
V.~K.~Beloshapka. Analytic complexity of functions of two variables.
\textit{Russian Journal of Mathematical Physics}, \textbf{14}:3 (2007), 243--249.
\texttt{doi:10.1134/S1061920807030016}.

\bibitem{Bel2016}
V.~K.~Beloshapka. Analytic complexity of functions of several variables.
\textit{Matematicheskie Zametki}, \textbf{100}:6 (2016), 781--789 (in Russian).
\texttt{doi:10.4213/mzm11320}.

\bibitem{Step2019}
M.~A.~Stepanova. On the analytic complexity of differential-algebraic functions.
\textit{Matematicheskii Sbornik}, \textbf{210}:12 (2019), 120--135 (in Russian).
\texttt{doi:10.4213/sm9146}.
English translation: \textit{Sbornik: Mathematics}, \textbf{210}:12 (2019), 1774--1787.

\bibitem{AHS}
R.~Ait El Manssour, M.~H\"ark\"onen, B.~Sturmfels. Linear PDE with constant coefficients.
\textit{Glasgow Mathematical Journal}, \textbf{65}, special issue S1 (2023), S2--S27.
\texttt{doi:10.1017/S0017089521000355}.

\bibitem{Robertz}
D.~Robertz. Recent progress in an algebraic analysis approach to linear systems.
\textit{Multidimensional Systems and Signal Processing}, \textbf{26}:2 (2015), 349--388.
\texttt{doi:10.1007/s11045-014-0280-9}.

\bibitem{Mal1956}
B.~Malgrange. Existence et approximation des solutions des \'equations aux d\'eriv\'ees
partielles et des \'equations de convolution.
\textit{Annales de l'Institut Fourier}, \textbf{6} (1956), 271--355.

\bibitem{Ehr1970}
L.~Ehrenpreis. \textit{Fourier Analysis in Several Complex Variables}.
Wiley-Interscience, New York, 1970.

\bibitem{Pal1970}
V.~P.~Palamodov. \textit{Linear Differential Operators with Constant Coefficients}.
Springer-Verlag, Berlin, 1970.

\bibitem{CHS2021}
Y.~Cid-Ruiz, R.~Homs, B.~Sturmfels. Primary ideals and their differential equations.
\textit{Foundations of Computational Mathematics}, \textbf{21}:5 (2021), 1363--1399.
\texttt{doi:10.1007/s10208-020-09485-6}.

\bibitem{AM1969}
M.~F.~Atiyah, I.~G.~Macdonald, Introduction to Commutative Algebra,
Addison--Wesley, 1969.

\bibitem{Bel2012}
V.~K.~Beloshapka. Analytical complexity: Development of the topic.
\textit{Russian Journal of Mathematical Physics}, \textbf{19}:4 (2012), 428--439.
\texttt{doi:10.1134/S1061920812040036}.

\bibitem{Bel2021}
V.~K.~Beloshapka. On simple solutions of the Burgers and Hopf equations.
\textit{Izvestiya Rossiiskoi Akademii Nauk. Seriya Matematicheskaya},
\textbf{85}:3 (2021), 5--12 (in Russian).
\texttt{doi:10.4213/im9051}.

\end{thebibliography}
\end{document}